\documentclass[leqno,11pt]{article}

\usepackage{amssymb,hyperref,amsthm}
\usepackage{euscript}
\usepackage[dvips]{graphicx}
\usepackage{flafter}
\usepackage{pstricks}
\usepackage[latin1]{inputenc}
\usepackage{amsmath}
\usepackage{amsfonts}
\usepackage{pstricks-add}
\usepackage{dsfont}
\usepackage{bm}
 \usepackage{lipsum}

\newcommand{\N}{{\mathbb N}}

\newcommand{\Q}{{\mathbb Q}}

\renewcommand{\d}{\text{ d}}
\newcommand{\dd}{\text{d}}

\def\ee{\varepsilon}

\def\L{{\bf L}}
\def\E{{\mathbb E}}

\def\P{{\mathbb P}}
\def\N{{\mathbb N}}

\usepackage[refpage]{nomencl}
\makenomenclature

\usepackage{makeidx}
\makeindex

\usepackage{mathrsfs} 

\evensidemargin\oddsidemargin

\usepackage{tikz}

\makeatletter \newcommand \listoftodos{\section*{Todo list} \@starttoc{tdo}}
\newcommand\l@todo[2]
{\par\noindent \textit{#2}, \parbox{10cm}{#1}\par} \makeatother

\begin{document}

\newcommand{\eqnsection}{
\renewcommand{\theequation}{\thesection.\arabic{equation}}
   \makeatletter
   \csname  @addtoreset\endcsname{equation}{section}
   \makeatother}
\eqnsection

\def\r{{\mathbb R}}
\def\e{{\mathbb E}}
\def\p{{\mathbb P}}
\def\P{{\bf P}}
\def\Pp{{ \mathbb{P}}}
\def\E{{\bf E}}
\def\Ee{{\mathbb E}}
\def\Q{{\bf Q}}
\def\z{{\mathbb Z}}
\def\N{{\mathbb N}}
\def\T{{\mathbb T}}
\def\G{{\mathbb G}}
\def\L{{\mathbb L}}
\def\1{{\mathds{1}}}
\def\deg{\chi}

\def\d{ \mathtt{d}}

\def\ss{ {\bf s} }
\def\d{\mathtt{d}}
\def\ttheta{{\bm \theta}}
\def\t{{\bf{t}}}
\def\a{{\bf{a}}}
\def\deg{\chi}
\def\B{\mathfrak{B}}

\def\M{{\mathbb{M}}}
\def\ee{\mathrm{e}}
\def\d{\, \mathrm{d}}
\def\S{\mathscr{S}}
\def\bs{{\tt bs}}
\def\bbeta{{\bm {S}}}
\def\ttheta{{\bm \theta}}

\newtheorem{theorem}{Theorem}[section]

\newtheorem{definition}[theorem]{Definition}
\newtheorem{Lemma}[theorem]{Lemma}
\newtheorem{Proposition}[theorem]{Proposition}
\newtheorem{Remark}[theorem]{Remark}
\newtheorem{corollary}[theorem]{Corollary}



\vglue50pt

\centerline{\Large\bf The Seneta-Heyde scaling for homogeneous fragmentations}

\bigskip
\bigskip


\medskip

\centerline{Andreas E. Kyprianou and Thomas Madaule}

\medskip

\centerline{\it University of Bath and Universit\'e Paul Sabatier}

\bigskip
\bigskip
\bigskip

{\leftskip=2truecm \rightskip=2truecm \baselineskip=15pt \small

\noindent{\slshape\bfseries Abstract.} Homogeneous mass fragmentation processes describe the evolution of a unit mass that breaks down randomly into pieces as time. Mathematically speaking, they can be thought of as continuous-time analogues of branching random walks with non-negative displacements. Following recent developments in the theory of branching random walks, in particular the work of \cite{AShi10}, we consider the problem of the Seneta-Heyde norming of the so-called additive martingale at criticality. Aside from replicating results for branching random walks in the new setting of fragmentation processes,  our main goal is to present a style of reasoning, based on $L^p$ estimates,  which works for a whole host of different branching-type processes. We show that our methods apply equally to the setting of branching random walks, branching Brownian motion as well as Gaussian multiplicative chaos. 

} 

\bigskip
\bigskip

\tableofcontents
\section{Introduction}

There exists a very strong mathematical analogy between fragmentation processes and  branching random walks; see for example \cite{bertfragbook, BRo04}. Both model a process of splitting. Where as, in the former, one understands `splitting' as the fragmentation of mass, in the latter, `splitting' corresponds to the creation of mass through offspring. Nonetheless, the genealogical tree-like structure that is inherently embedded in the stochastic evolution of both processes accounts for many mathematical similarities. This is not only the case for fragmentation processes and branching random walks. The right mathematical perspective also reveals commonly embedded structures in branching particle diffusions, superprocesses and multiplicative chaos. 

An object that can be found commonly amongst all of the aforesaid processed is the so-called additive martingale. Roughly speaking, representing the configuration of any of the aforementioned processes at a fixed time as a measure valued object (atomic measures in the case of fragmentation processes, branching random walks and branching Brownian motion),  the additive martingale emerges by integrating the aforesaid measure against an exponential function in space and exponentially discounting in an appropriate way in time. The additive martingale packages information about the empirical `spatial configuration' of the process in question into a convenient stochastic process whose limit can be guaranteed in the almost sure sense (thanks to positivity). Traditionally, the martingale limit, when non-trivial, has played a very important role in studies pertaining to the asymptotic `spatial' evolution of such branching-type processes.

The main objective of this paper is to consider the problem of the Seneta-Heyde norming of the additive martingale for a critical parameter regime in the setting of fragmentation processes. That is to say, for a particular extreme choice of parameters, the additive martingale will tend almost surely to zero (thereby offering little direct information about the process in the limit) whereupon we establish a deterministic renormalisation of the martingale such that it has a non-trivial limit. Moreover, in doing so, we will present a style of reasoning, based on $L^p$ estimates,  which works for a whole host of different branching-type processes. We show that our methods apply to the setting of branching random walks, branching Brownian motion as well as Gaussian multiplicative chaos. 

The rest of the paper is outlined as follows. In the next section we introduce fragmentation processes. In section \ref{mainresultssection} we present our main results. The section thereafter reviews some technical considerations that will be used in Section \ref{sectionproofs}, where the main results are proved. In Section \ref{BRWBBMMC} we briefly review how the methodology transfers to the case of branching random walks, branching Brownian motion and Gaussian multiplicative chaos. Finally, in the Appendix, we list some results for L\'evy processes which are needed in the proofs of the main results, but appear not to be contained in existing literature, despite their analogues for random walks being available. 


\section{Notations and definitions}
Below we give a brief overview of the definition and structure of a homogeneous fragmentation process. The reader is referred to Bertoin \cite{bertfragbook} for more details.
Let $\mathcal{P}$ be the space of partitions of the natural numbers. Here a partition of $\mathbb{N}$ is a sequence $\pi=(\pi_1, \pi_2, \cdots)$ of disjoint sets, called blocks, such that $\bigcup_{i\in\N} \pi_i = \mathbb{N}$. The blocks of a partition are enumerated in the increasing order of their least element, that is to say $\min \pi_i\leq \min\pi_j$ when $i\leq j$ (with the convention that $\min \emptyset = \infty$). 
Now consider the measure $\mu$ on $\mathcal{P}$, given by
\begin{equation}
\mu(d\pi) = \int_{\mathcal{S}}\varrho_{\bf s}(d\pi)\nu(d{\bf s}),
\label{implicit}
\end{equation}
where $\varrho_{\bf s}$  is the law of Kingman's paint-box based on ${\bf s}\in\mathcal{S}$ (cf. page~98 of Bertoin \cite{bertfragbook}) with
 \[
\mathcal{S}: = \left\{\mathbf{s}=(s_1, s_2, \cdots) : s_1\geq  s_2 \geq \cdots \geq 0 , \, \sum_{i=1}^\infty s_i \le 1\right\},
 \] 
and the so-called {\it dislocation measure}, $\nu\neq 0$,  is a measure on $\mathcal{S}$
 satisfying
\begin{equation} 
\nu({\bf s}\in\mathcal{S}:s_2=0)=0
\label{e.levymeasure.0}
\end{equation}
as well as
\begin{equation}
\int_{\mathcal{S}}(1-s_1)\nu(d{\bf s})<\infty.
\label{e.levymeasure}
\end{equation}
 It is known that $\mu$ is  an exchangeable partition measure, meaning that it is invariant under the action of finite permutations on $\mathcal{P}$. It is also known (cf. Chapter 3 of Bertoin \cite{bertfragbook}) that it is possible to construct a fragmentation process  on the space of partitions $\mathcal{P}$ with the help of  a Poisson point process $\{(\pi(t), k(t)) : t\geq 0\}$ on $\mathcal{P}\times \mathbb{N}$ which has intensity measure $\mu\otimes\sharp$, where $\sharp$ is the counting measure. The aforementioned   $\mathcal{P}$-valued fragmentation process is a  Markov process which we denote by  $\Pi = (\Pi(t))_{  t\geq 0}$, where
 $\Pi(t) = (\Pi_1(t), \Pi_2(t),\cdots)\in\mathcal{P}$ 
 is such that at all times $t\geq0$ for which an atom $(\pi(t), k(t))$ occurs in $(\mathcal{P}\backslash(\mathbb{N},\emptyset,\ldots))\times \mathbb{N}$, $\Pi(t)$ is obtained from $\Pi(t-)$ by partitioning the $k(t)$-th block  into the sub-blocks $(\Pi_{k(t)}(t-) \cap \pi_j(t) : j=1,2,\cdots)$. When $\nu$ is a finite measure each block experiences an exponential holding time before it fragments. 
 For future reference, we write  $\mathbb{G}: = (\mathcal{G}(t))_{t\geq 0}$ for the natural filtration associated to the Poisson point process generating the fragmentation process.

Thanks to the properties of the exchangeable partition measure $\mu$ it can be shown that, for each $t\geq 0$, the distribution of $\Pi(t)$ is exchangeable and that the  blocks of $\Pi(t)$  have asymptotic frequencies in the sense that, for each $i\in\mathbb{N}$, the limit 
\[
|\Pi_i(t)|:=\lim_{n\to\infty} \frac{1}{n}\sharp\{\Pi_i(t)\cap \{1,\cdots, n\} \}
\]
exists almost surely.  Moreover, it is known that $|\Pi_i(t)|$ exists $\mathbb P$--a.s. simultaneously for all $t\ge0$ and $i\in\mathbb N$. 
We can thus define the {\it mass fragmentation process} $X(t) = (X_1(t), X_2(t), \cdots)$, $t\geq 0$, where $X_i(t)$ is the mass of the $i$-th largest block in the sequence $\Pi(t)$, when blocks  are re-ordered by their asymptotic frequencies.  

The mass fragmentation process $X(t)$, $t\geq 0$, can equivalently be constructed using the projection  to the space $\mathcal{S}\times \mathbb{N}$ of the Poisson point process that generates $\Pi$. This is a Poisson point process with intensity $\nu\otimes\sharp$ where  $\nu$ is related to $\mu$ via the relation (\ref{implicit}). 
Suppose that $({\bf s}(t), k(t))$ is an atom  in $\mathcal{S}\backslash\{(1,0,\cdots)\}\times \mathbb{N}$ in this Poisson point process that occurs at time $t\geq 0$.  The sequence ${X}(t)$ is obtained from ${X}(t-)$ by replacing its $k(t)$-th term, $X_{k(t)}(t-)$, with the sequence $X_{k(t)}(t-){\bf s}(t)$ and ranking all existing blocks and new blocks in decreasing order.  The process $(X(t))_{t\geq0}$ is now adapted to the subfiltration of $\mathbb{G}$, say   $\mathbb{F}: = (\mathcal{F}(t))_{t\geq 0}$, which is generated by the Poisson point process on $\mathcal{S}\times \mathbb{N}$.

As a consequence of either of the above two descriptions, the mass fragmentation process has a convenient Markov property (also known as the {\it fragmentation property}) as follows. Given that $X(t)  = (s_1,s_2,\cdots)\in\mathcal{S}$, where $t\geq 0$,  for $u>0$, $X(t+u)$ has the same law as the process obtained by ranking in decreasing order the sequences $X^{(1)}(u), X^{(2)}(u),\cdots$ where the latter are independent, random mass partitions with values in $\mathcal{S}$ having the same distribution as $X(u)$ but scaled in size by the factors under $s_1$, $s_2$, $\cdots$ respectively.

Let us define $\xi_t  := -\log |\Pi_1(t)| = -\log X_1(t)$,  $t\geq 0$, with the convention $-\log0:=\infty$. The process $(\xi_t)_{t\geq 0}$ is called the {\it tagged fragment}. Using the Poissonian construction of the fragmentation process, one easily shows that $(\xi_t)_ { t\geq 0}$ is a killed
subordinator with cemetery state $\infty$ and killing rate
\begin{equation}\label{e.kappa}
\kappa:=\int_{\mathcal S}\left(1-\sum_{k\in\N}s_k\right)\nu(\dd{\bf s})\geq 0.
\end{equation}
Moreover, it is well known that its Laplace exponent $\Phi$, given by 
\[
e^{-\Phi(p)t} : = \mathbb{E}(e^{-p\xi_t}),\qquad t\geq 0,
\] 
is well defined over $[\underline{p},\infty)$, where 
\[
\underline{p}:= \inf\left\{ p\in \r: \int_{\mathcal{S}} \sum_{i=2}^\infty s_i^{p+1}\nu(d \ss)<\infty\right\}\in[-1,0]
\]
For convenience, we assume that 
\begin{equation}
\underline{p}>-1.
\label{Asump2}
\end{equation}
Then, for every $q>\underline{p}$,
\begin{equation}
\Phi(q)= \int_{\mathcal{S}}(1-\sum_{i=1}^\infty s_i^{q+1})\nu(d\ss).
\end{equation}
The function $\Phi$ is a concave infinitely differentiable  increasing function. We shall assume that the equation 
\begin{equation}
\Phi(q)=(q+1)\Phi'(q)
\text{ admits a unique solution } \bar{p}>\underline{p}
\label{Asump1}
\end{equation}
 Such an assumption is automatically satisfied if there exists some $p^*$ such that $\Phi(p^*)=0$;  see \cite{Ber03}. In particular this assumption is automatically satisfied in the conservative case (i.e. $\nu(\sum_{i\geq 1} s_i <1)=0$). More precisely, we have
\begin{equation}
\Phi(q)-(q+1)\Phi'(q)<0\Longleftrightarrow q\in ]\underline{p},\bar{p}[
\end{equation}
and 
\begin{equation}
\text{the map } q\to \Phi(q)/(q+1) \text{ increases on }]\underline{p},\bar{p}[ \text{ and decreases on } ]\bar{p}, \infty[.
\end{equation}

The main object of interest in this paper is the following family of random processes:
\begin{equation}\label{addmg}
W(t,p):= \ee^{\Phi(p)t}\sum_{i=1}^\infty X_i(t)^{1+p},\qquad t\geq 0,
\end{equation}
with $ p>\underline{p}$.
From Bertoin and Rouault \cite{BRo04}, we know that $W(t,p)$ is a positive martingale (and hence almost surely convergent with limit, say $W(\infty,p)$) having càdlàg paths and mean $1$. In particular for any $p>\underline{p}$ and $t\geq 0$,
\begin{equation}
\label{relationPhi}
\E(\sum_{i=1}^\infty X_i(t)^{1+p})= \exp(-t\Phi(p))<\infty.
\end{equation} 
Moreover when $t$ goes to infinity, we know that $ W(t,p)\to W(\infty,p)>0$ almost surely if $ p\in (\underline{p},\bar{p})$ and $ W(t,p)\to 0$ almost surely if $p\geq \bar{p} $. Here we are interested, in the critical case, i.e in the behaviour of $W(t,\bar{p})$. Our study involves the so-called {\it critical derivative martingale}
\begin{equation}\label{daddmg}
M'(t):=  \ee^{\Phi(\bar{p}) t} \sum_{i=1} \{- \log X_i(t)-t\Phi'(\bar{p}) \} X_i^{1+\bar{p}}(t), \qquad t\geq 0,
\end{equation}
which is of interest because it converges almost surely to a non-trivial value. Indeed, 
 Bertoin and Rouault \cite{BRo04} and the preceding unpublished preprint \cite{BRpreprint} show that 
\begin{equation}
M'(\infty): = \lim_{t\to\infty } M'(t)
\end{equation}
exists almost surely with $M'(\infty)$ a strictly  positive random variable.

\section{Main results}\label{mainresultssection}

Our first main result is an analogue to the so-called Seneta-Heyde normalization at criticality for branching random walks. The result in that setting is due to \cite{AShi11} and  shows that the additive martingale tends to zero at the exact rate of $t^{-1/2}$. More precisely, it shows that, up to a multiplicative constant, when multiplying the additive martingale by the square root of time, there is weak convergence to the derivative martingale limit. 

\begin{theorem}
\label{mainresult}Let $\sigma^2:=\Phi''(\bar{p})$, which is finite and positive thanks to the smoothness and concavity of $\Phi$. 
Assuming (\ref{Asump1}) and (\ref{Asump2}), we have 
\begin{equation}
\label{seneta}
\underset{t\to \infty}{\lim} \sqrt{t}W(t,\bar{p})= \sqrt{\frac{2}{\pi\sigma^2}}M'( \infty),\quad \text{in } \P \text{-probability.}
\end{equation}
\end{theorem}

\bigskip

In light of the conclusion above, as well as  the expressions for the additive and derivative martingales given in  (\ref{addmg}) and (\ref{daddmg}) respectively, one might be tempted to compare the terms $- \log X_i(ut) -ut\Phi'(\bar{p})$, $i\in\mathbb{N}$, in the derivative martingale, with $\sqrt{t}$. To this end,
let $\mathcal{D}([0,a])$, for $a>0$, be the space of càdlàg functions over the time horizon $[0,a]$ which we equip with the Skorokhod topology. Observe that for any $t>0$, $i\in \N$, 
\begin{equation}
{t^{-1/2}}(- \log X_i(ut) -ut\Phi'(\bar{p}))_{u\in [0,1]}\in \mathcal{D}([0,1]).
\end{equation}
Write 
 $ \mathcal{C}(\mathcal{D},[0,a])$ and $\mathcal{C}_b(\mathcal{D},[0,a])$ 
 for the space of non-negative, continuous (with respect to the Skorokhod topology) functions on $\mathcal{D}([0,a])$ and the space of bounded functions in $\mathcal{C}(\mathcal{D},[0,a])$ respectively.
 
 For convenience, let us also introduce the notation 
 \[
  Z_i(t)= -\log(X_i(t)) -t\Phi'(\bar{p}) = -\log(X_i(t)) -t\frac{\Phi(\bar{p})}{1+\bar{p}}  \qquad t\geq 0
 \]
 where the second equality is a consequence of the definition of $\bar{p}$ and the index $i\in\mathbb{N}$ is the index set of fragments at time $t$.
\begin{theorem}
\label{resultBES}
Let $(R_s)_{s\in [0,1]}$ a three dimensional Bessel process. Assuming (\ref{Asump1}) and (\ref{Asump2}), for any functional $F\in \mathcal{C}_b(\mathcal{D},[0,1])$, we have
\begin{equation}
\lim_{t\to \infty} \sum_{i=1} {Z_i(t)}\ee^{- (1+\bar{p}) Z_i(t)}F( {t^{-1/2}}({Z_i(ut)})_{u\in [0,1]}  )   = M'(\infty) \E\left(F((\sigma R_s)_{s\in [0,1]}) \right),\qquad\text{in } \P\text{-probability}.\label{eqresultBES} 
\end{equation}
\end{theorem}

\bigskip

Combining the last two theorems, one essentially has a functional law of large numbers. Specifically, for any $F\in \mathcal{C}_b(\mathcal{D},[0,1])$,
\[
\lim_{t\to \infty}  \sum_{i=1}  \left( \frac{\ee^{- (1+\bar{p}) Z_i(t)} }{W(t,\bar{p})}\right){t^{-1/2}} {Z_i(t)}F( {t^{-1/2}}({Z_i(ut)})_{u\in [0,1]}  )    =  \sqrt{\frac{\pi\sigma^2}{2}}\E\left(F((\sigma R_s)_{s\in [0,1]}) \right)
\]
in probability.

\section{Technical considerations for fragmentation processes}

In this section, we remind the reader of some straightforward standard theory for fragmentation processes that will be used in the later text. 
\subsection{Change of measures}

%

Recalling that $-\log|\Pi_1(t)|$, $t\geq0$ is a subordinator, we can appeal to the well-known fact for subordinators is that
\begin{equation*}
\mathcal{E}(t,p):= \ee^{\Phi(p)t}|\Pi_1(t)|^p,\qquad p>\underline{p},
\end{equation*}
is a positive $(\P,\mathbb{G})-$martingale. Moreover, when we project $\mathcal{E}(p,t)$ on the sub-filtration $\mathbb{F}$, we recover the additive martingale $ W(t,p)= \exp(t\Phi(p))\sum_{i=1}^\infty |\Pi_i(t)|^{p+1}= \exp(t\Phi(p))\sum_{j=1}^\infty X_j^{p+1}.$

Following an original methodology  of Lyons, Pemantle and Peres, but adapted to the setting of fragmentation process,  (see e.g \cite{Lyo97, BRo04}), we introduce the {\it tilted probability measure} $\P^{(p)}$ given by 
\begin{equation}
d\P^{(p)}_{\big| \mathcal{G}(t)}=\mathcal{E}(t,p)\,d\P_{\big| \mathcal{G}(t)},\qquad t\geq 0.
\end{equation}
Observe that projections on the sub-filtration $\mathbb{F}$ give the identity
\begin{equation}
d\P^{(p)}_{\big| \mathcal{F}(t)}= W(t,p)\, d\P_{\big| \mathcal{F}(t)}, \qquad t\geq 0.
\end{equation}
The effect of the change of probability is straightforward to describe, both at the level of the tagged fragment and that of the Poisson point process, say $N(dt, d\pi, k)$, on $[0,\infty)\times\mathcal{P}\times\mathbb{N}$ that generates the partition-valued fragmentation process.

\begin{Proposition}[Bertoin and Rouault \cite{BRo04}]
\label{BerRouau}  Fix $p>\underline{p}$.
Let $\mathcal{D}_1\subset [0,\infty[$ be the random set of times $r\geq 0$ for which the Poisson point process, $N(dt, d\pi, k)$, has an atom on the fiber $\{r\}\times \mathcal{P}\times \{1\}$, and for every $r\in \mathcal{D}_1$, denote the second component of this atom by $\pi(r)$.
\begin{itemize}
\item[(i)] Under $\P^{(p)}$, the process $\xi_t= -\log |\Pi_1(t)|$ is a subordinator with Laplace exponent 
$$ 
\Phi^{(p)}(q):= \Phi(p+q)-\Phi(p),\qquad q>\underline{p}-p.
$$
\item[(ii)] Under $\P^{(p)}$: 
\begin{itemize}
\item The restriction of $N(dt, d\pi, k)$ to $\r_+\times \mathcal{P}\times \{2,3,...\}$ has the same law as it has under $\P$ and is independent of the restriction to the fiber $\r_+\times \mathcal{P}\times \{1\}$.

\item The restriction of $N(dt, d\pi, k)$ to $\r_+\times \mathcal{P}\times \{1\}$ is a Poisson point process on  $\r_+\times \mathcal{P}$ with intensity $dr\otimes \mu^{(p)}$, where 
$$
\mu^{(p)}(d\pi)= |\pi_1|^p\mu(d\pi).
$$

\end{itemize}
\end{itemize}
\end{Proposition}
From Proposition \ref{BerRouau}, we deduce the following corollary for the particular case that $p = \bar{p}$.

\begin{corollary}[Many-to-one Lemma]
	\label{Manytoone}
For any  $F\in \mathcal{C}(\mathcal{D},[0,t])$ and $p>\underline{p}$, we have
\begin{equation}
\E\left(\sum_{i=1}^\infty\ee^{\Phi({p})t}  X_i(t)^{1+\bar{p}} F\Big([X_i(st)]_{s\in [0,1]}\Big)  \right)=\E^{({p})}\left[ F( [{\rm e}^{-\xi_{st}}]_{s\in [0,1]})\right],\qquad t\geq 0.
\end{equation} 
\end{corollary}

\begin{Remark}\rm 
	\label{laremarke}
Another `many-to-one' type interpretation of the change of measure (previously observed in \cite{BRo04}, see also \cite{Lyo97})
is that, given $\mathcal{F}(t)$, we may choose a `tagged fragment' in an empirical way with probabilities proportional to $\ee^{\Phi({p})t} X_i(t)^{1+{p}}$. Moreover, the historical evolution of the tagged fragment satisfies
 \[
 F([{\rm e}^{-\xi^*_{st}}]_{s\in[0,1]}) := \sum_{i=1}^\infty \frac{\ee^{\Phi({p})t} X_i(t)^{1+{p}} }{W(t, p)}F\Big([X_i(st)]_{s\in [0,1]}\Big), \qquad t\geq 0,
\]
for   $F\in \mathcal{C}(\mathcal{D},[0,t])$. Moreover, under $\mathbf{P}^{({p})}$, $(\xi^*)_{t\geq 0}$ is a subordinator with Laplace exponent $\Phi^{(p)}$.

\end{Remark}

\begin{Remark} \rm
Observe that $\Phi'(\bar{p})= \E^{(\bar{p})}(\xi_1)$, and that  $\sigma^2$ can be rewrite as
\begin{equation}
\E^{(\bar{p})}\left( [\xi_1-\E^{(\bar{p})}(\xi_1)]^2\right) = \E\left( \sum_{i=1}^\infty \{\log X_i(1)+ \Phi'(\bar{p})\}^2 X_i(1)^{1+\bar{p}}\ee^{\Phi(\bar{p})}\right)<\infty.
\end{equation}
Moreover observe that, for any $\epsilon>0$ such that  $\epsilon <\bar{p}-\underline{p}$,
\begin{equation}
\E^{(\bar{p})}(\ee^{\epsilon \xi_t})=\E\left(\ee^{\Phi(\bar{p})t} \sum_{i=1}^\infty X_i(t)^{1+\bar{p}} X_i(t)^{-\epsilon}  \right)<\infty, \qquad t\geq 0.
\end{equation}
Hence, the subordinator $(\xi_t)_{t\geq 0}$ under $\P^{(\bar{p})}$ has moment positive and negative exponential moments. Indeed, for any $t>0$ and $\epsilon <\bar{p}-\underline{p}$ , we have
\begin{equation}
\label{momentexp}
\E^{(\bar{p})}\left( \ee^{\epsilon |\xi_t|}  \right)<\infty.
\end{equation}
\end{Remark}

\subsection{On the integrability of the dislocation measure $\nu$}
In the following Lemma we point out an integrability property for the fragmentation process which will be essential in the proof of  Lemma \ref{Momentdomination}. 
\begin{Lemma}
	\label{LEintegrab}
	Under (\ref{e.kappa}), (\ref{Asump2}) and (\ref{Asump1}), for any $q>\underline{p}>-1$ and $\gamma\in (1,2)$ such that $\gamma q>  \underline{p}>-1$ we have
	\begin{equation}
		\label{integrab} \E\Big( \Big(  \sum_{i=1}^{\infty} X_i(t)^{1+q}\Big)^\gamma  \Big) <\infty,\qquad \forall t\geq 0.
	\end{equation}
\end{Lemma}
\noindent{\it Proof of Lemma \ref{LEintegrab}.} Fix $t>0$. By (\ref{relationPhi}) for any $p>-1$, the inequalities
$$ \E\Big(\sum_{i=1}^\infty X_i(t)^{1+p}\Big)<+\infty,\qquad \int_{\mathcal{S}}  \sum_{i=2}^\infty s_i^{1+p} \nu(d {\bf s}) <\infty $$
are both equivalent.  As the function $p\mapsto \sum_{i=1}^\infty {X_i(t)^p}$ is decreasing and  $\gamma q >\underline{p}$, we thus get
\begin{equation}
	\label{musumE}\E\left(    \left(  \sum_{i=1}^{\infty} X_i(t)^{ {1+\gamma q}   }       \right)\right)<+\infty.
\end{equation}  
Now it suffices to apply the Hölder inequality with $a=\gamma$ and $b={\gamma \over \gamma-1}$, to get
\begin{eqnarray*}
	\E\left( \left(  \sum_{i=1}^{\infty} X_i(t)^{1+q}\right)^\gamma  \right)&= &\E\left(    \left(  \sum_{i=1}^{\infty} X_i(t)^{\frac{1+\gamma q}{\gamma}  }    X_i(t)^{\frac{\gamma-1}{\gamma}}    \right)^\gamma       \right)
	\\
	&\leq & \E\left(    \left(  \sum_{i=1}^{\infty} X_i(t)^{ {1+\gamma q}   }       \right) \left(  \sum_{i=1}^{\infty}     X_i(t)    \right)^{\gamma-1}      \right) \\
	&\leq& \E\left(    \left(  \sum_{i=1}^{\infty} X_i(t)^{ {1+\gamma q}   }       \right)\right)<+\infty,
\end{eqnarray*}
where we used $\sum_{i=1}^\infty X_i(t)\leq 1$ in the second inequality. 
\hfill$\Box$

\begin{Remark}\rm By exactly the same arguments we also have
\begin{eqnarray}
	\nonumber   \int_{\mathcal{S}}   \left(  \sum_{i=2}^{\infty} s_i^{1+q}\right)^\gamma  \nu(d {\bf s}) &= & \int_{\mathcal{S}}    \left[ \left(  \sum_{i=2}^{\infty} s_i^{( 1+\gamma q)/{\gamma}}  s_i^{(\gamma-1)/{\gamma}}\right)^\gamma\right] \nu(d {\bf s})
	\\
	\label{holder2}&\leq&  \int_{\mathcal{S}}     \left( \sum_{i=2}^\infty s_i^{  1+\gamma q}\right) \left( \sum_{i=1}^\infty s_i \right)^{\gamma-1} \nu(d {\bf s})\notag\\
	&\leq&       \int_{\mathcal{S}}   \sum_{i=2}^\infty s_i^{ 1+\gamma q}  \nu(d {\bf s})< +\infty
\end{eqnarray}
\end{Remark}

\subsection{The truncated processes}
Studying $(W(t,\bar{p}))_{t\geq 0}$ or $(M'(t))_{t\geq 0}$ is hard because both processes are not well concentrated around their mean. To overcome this difficulty we will study instead the truncated versions of these martingales. The idea of truncating is now standard, and has successfully been used in the setting of branching random walks and branching Brownian motion, e.g.  \cite{Harris-KPP, Kyp-KPP, BK04}, as well as in the current setting, e.g. \cite{BRo04}. Recall that, for $t\geq 0$, $Z_i(t): = -\log X_i(t) - \Phi'(\bar{p})t$, where $i\in\mathbb{N}$ is the index set of blocks at time  $t$. For any $a>0$ and $i\in \N$, let $\zeta^a_i:= \inf\{t\geq 0,\, Z_i(t)<-a\}$. For any $t\in \r^+$ and $F\in \mathcal{C}(\mathcal{D},\r_+)$ we define:
\[
	W^{(a)}(t,\bar{p}):= \sum_{i=1}^{\infty} \ee^{- (1+\bar{p}) Z_i(t)}\1_{\{{\zeta^a_i}> t \}},  \qquad
	M'^{(a)}(t):=\sum_{i=1} Z_i(t)  \ee^{- (1+\bar{p}) Z_i(t)}\1_{\{{\zeta^a_i}> t \}},
	\]
	and 
	\[
	M'^{(a)}(t,F):=  \sum_{i=1} (a+ Z_i(t))  \ee^{- (1+\bar{p}) Z_i(t)}\1_{\{{\zeta^a_i}> t \}}
	 F\{ {t}^{-1/2}(Z_i(ut))_{u\in [0,1]} \}.
	 \]
Notice that when $a$ goes to infinity we recover the originally processes. Indeed this follows from the result in \cite{BRo04} that  $\inf_{i\in \N}\inf_{t\geq 0}( -\log(X_i(t)) -t\Phi'(\bar{p}) )>-\infty$, a.s. In particular,  it was show in the aforesaid paper that
\begin{equation}
	\label{convDeriv}
	\lim_{t\to\infty} M'^{(a)}(t)=M'^{(a)}(\infty)>0,\qquad \text{a.s}.
\end{equation}
%
%

\section{Proofs of Theorem \ref{mainresult} and \ref{resultBES}}\label{sectionproofs}
The proofs of the Theorem \ref{mainresult} and \ref{resultBES} are very similar. They are combinations of the two followings propositions.
\begin{Proposition}
	\label{ConditionalE}Fix $a>0$.
	Assume  (\ref{e.kappa}), (\ref{Asump2}) and (\ref{Asump1}).
	\begin{itemize}
	\item[(i)] For any $l>0$,  we have
	\begin{equation}
		\label{eqConditionali}  \underset{t\to\infty}{\lim} \sqrt{t}\, \E\left( W^{(a)}(t+l,\bar{p})\big|\mathcal{F}_l\right) \overset{\text{a.s}}{=} \sqrt{\frac{2}{\pi\sigma^2}} M'^{(a)}({l}).
	\end{equation}
	
	\item[(ii)]Let $(R_s)_{s\geq 0}$ a three dimensional Bessel process. For any $F\in \mathcal{C}(\mathcal{D},\r^+)$ and $l\geq 0$, we have
	\begin{equation}
		\label{eqConditionalii}  \underset{t\to\infty}{\lim} \E\left( M'^{(a)}(t+l,F)\big|\mathcal{F}_l\right) \overset{\text{a.s}}{=} \E\left(F((\sigma R_s)_{s\in [0,1]}) \right)   M'^{(a)}(l) .
	\end{equation}
	\end{itemize}
\end{Proposition}

\begin{Proposition}
	\label{Concentration}Fix $a>0$.
	Assume  (\ref{e.kappa}), (\ref{Asump2}) and (\ref{Asump1}). There exists $c>0$, $\gamma\in (1,2)$ such that:
	
	\begin{itemize}
	\item[(i)] For any $\epsilon >0,\, l,\, t\geq 0$ we have
	\begin{equation}
		\label{SeneCond}
		\P\left(  \Big| \E\left( (t+l)^{1/2}\, W^{(a)}(t+l,\bar{p}) \big|\mathcal{F}_l\right)-   (t+l)^{1/2}\, W^{(a)}(t+l,\bar{p})\Big|\geq \epsilon \right)\leq \frac{c(t+l)^\frac{\gamma}{2}}{\epsilon^\gamma t^{\frac{\gamma}{2}} \sqrt{l}}(1+a) \ee^{(\gamma-1) a}.
	\end{equation}
	
	\item[(ii)]For any $\epsilon >0,\, l,\, t\geq 0$ and $F\in \mathcal{C}(\mathcal{D},\r_+)$ we have
	\begin{equation}
		\label{BesCond}
		\P\left(  \Big| \E\left( M'^{(a)}({t+l}, F) \big|\mathcal{F}_l\right)-   M'^{(a)}(t+l,F)\Big|\geq \epsilon \right)\leq \frac{c}{\epsilon^\gamma \sqrt{l}}(1+a) \ee^{(\gamma-1) a}.
	\end{equation}
\end{itemize}
\end{Proposition}
The proof of the Proposition \ref{ConditionalE} only involves the many-to-one Lemma (Corollary \ref{Manytoone}). Proposition \ref{Concentration} is a concentration result, for which the integrability property (\ref{integrab}) is essential for its proof. 
\\

\noindent{\it Proofs of the theorems \ref{mainresult} and \ref{resultBES}.}
By combining \eqref{eqConditionali}, \eqref{SeneCond} and (\ref{convDeriv}) we have
\begin{equation}
	\label{combine1} \lim_{t\to\infty} \sqrt{t} W^{(a)}(t)= \sqrt{\frac{2}{\pi\sigma^2}} M'^{(a)}(\infty),\qquad \text{in } \P \text{ probability.}
\end{equation}
As alluded to earlier, from \cite{BRo04}, their result that  $\inf_{i\in \N}\inf_{t\geq 0}( -\log(X_i(t)) -t\Phi'(\bar{p}) )>-\infty$ a.s. implies that
\begin{equation}
	\label{combine2} \lim_{a\to\infty } M'^{(a)}(\infty)\overset{a.s}{=} M'(\infty)\qquad\text{and    }\quad \lim_{a\to\infty } \lim_{t\to\infty} \left| \sqrt{t} W^{(a)}(t)-\sqrt{t}W(t)\right|= 0.
\end{equation}
Finally Theorem \ref{mainresult} is  a direct consequence of  (\ref{combine1}) and (\ref{combine2}). The proof of Theorem \ref{resultBES} is identical, it suffices to use (\ref{eqConditionalii}) and (\ref{BesCond}), instead of \eqref{eqConditionali} and \eqref{SeneCond}.\hfill$\Box$

\subsection{Computation of the conditional expectation}
This sub-section is dedicated to the proof of the Proposition \ref{ConditionalE}. Throughout, we appeal to results from the Appendix, all of which concern asymptotic distributional properties of spectrally one-sided  L\'evy processes. The reader is encouraged to briefly browse the results there before reading on, all of which are known for random walks but are seemingly missing from the L\'evy process literature, thereby necessitating a proof.   We start by proving \eqref{eqConditionali}.
\\

\noindent{\it Proof of (\ref{eqConditionali}).} Let us define the process $(Z_t)_{t\geq0}$ with probabilities  $\mathbb{P}_x$, $x\in\mathbb{R}$ such that, under $\mathbb{P}_x$, it is equal in law to 
$ \xi_t - t\Phi'(\bar{p})$, $t\geq 0$, under $\P^{(\bar{p})}$.
Note that it is a spectrally positive L\'evy process of bounded variation with mean  $0$ and variance $\sigma^2$. For convenience, we shall also write 
\[
\zeta^a = \inf\{t>0 : Z_t < -a\}, \qquad a\in \mathbb{R}.
\]
According to the fragmentation property, then Corollary \ref{Manytoone}, for any $t,l>0$, 
\begin{eqnarray*}
	&&\E\left( \sum_{i=1}^{\infty} \ee^{-(1+\bar{p})Z_i(t+l)}\1_{\{{\zeta^a_i}> t+l \}}  \big|\mathcal{F}_l\right)
	\\
	&&=  \sum_{i=1}^\infty \ee^{-(1+\bar{p})Z_i(l)}\1_{\{{\zeta^a_i}>l \}}  \E\left( \sum_{j=1}^\infty \ee^{-(1+\bar{p})Z_j(t)} \1_{\{\zeta^{a+x}_j>t\}}  \right)_{\big| x={Z_i(l)}}
	\\
	&&=  \sum_{i=1}^{\infty} \ee^{-(1+\bar{p})Z_i(l)} \1_{\{{\zeta^a_i}>l \}}  \P^{(\bar{p})}\left( \zeta^{a+x}>t \right)_{\big| x={Z_i(l)} }.
\end{eqnarray*}
Moreover by (\ref{Tool1}) we have for any $i\in \N$,
\begin{eqnarray*}
	&&\lim_{t\to\infty} \sqrt{t}\P^{(\bar{p})}\left( \zeta^{a+x} >t \right)_{\big| x={Z_i(l)}}=  \sqrt{\frac{2}{\pi\sigma^2}}(a+ {Z_i(l)} ).
\end{eqnarray*}
Hence,  we deduce that
\[
	\lim_{t\to\infty} \sqrt{t} \E\left( \sum_{i=1}^{\infty} \ee^{-(1+\bar{p})Z_i(t+l)}\1_{\{{\zeta^a_i}> t+l \}}  \big|\mathcal{F}_l\right)
	 =\sqrt{\frac{2}{\pi\sigma^2}} \sum_{i=1}(a+Z_i(l))  \ee^{-(1+\bar{p})Z_i(l)}\1_{\{{\zeta^a_i}>l \}} =\sqrt{\frac{2}{\pi\sigma^2}} M'^{(a)}(l),
\]
which concludes the proof of (\ref{eqConditionali}). \hfill$\Box$
\\

\noindent{\it Proof of (\ref{eqConditionalii}).} The proof works in a similar way to the proof of (\ref{eqConditionali}) but instead of (\ref{Tool1}) we will use Proposition \ref{Lemapart}. For any $t,l>0$, 
\begin{eqnarray}
	\nonumber  &&\E\left(  \sum_{i=1}^\infty (a+Z_i(t+l))\ee^{-(1+\bar{p})Z_i(t+l)}\1_{\{{\zeta^a_i}> t+l \}} F\{(t+l)^{-1/2}(Z_i(u(t+l))_{u\in [0,1]} \} \big|\mathcal{F}_l\right)
	\\
	\nonumber &&=  \sum_{i=1} \ee^{-(1+\bar{p})Z_i(l)}  \1_{\{{\zeta^a_i}>l \}}  \mathbb{E}\left(   ({Z_t}+a+g_l )\1_{\{ \zeta^{a+g_l}>t  \}}\times\right.
	\\
	\label{expectassion} &&  \hspace{5cm} F\left\{ { (t+l)^{-1/2}}\Big( g_{st}\1_{\{ s\leq \frac{l}{t} \}}   +Z_{st}\1_{\{ s\geq \frac{l}{t} \}}\Big)_{s\in [0,1]} \right\} \Big)_{\big| (g_s)_{s\leq l}=(Z_i(s))_{s\leq l} }.
\end{eqnarray}
Let us study the expectation inside the sum. 
\begin{eqnarray}
	&& \mathbb{E}_{a+g_l}\left( \frac{{Z}_t}{\sqrt{t}} F\left\{  { (t+l)^{-1/2}}\Big( g_{st}\1_{\{ s\leq \frac{l}{t} \}}  +{Z}_{st}    \1_{\{ s\geq \frac{l}{t} \}}\Big)_{s\in [0,1]} \right\}    \big|  \inf_{s\leq t} {Z}_s\geq 0  \right)\times {\sqrt{t}} \Pp_{a+g_l}( \inf_{s\leq t} {Z}_s\geq 0)
\end{eqnarray}            
By (\ref{Tool1}), we know that $\lim_{t\to \infty} \frac{ (t+l)^{1/2}}{a+g_l} \Pp_{a+g_l}( \inf_{s\leq t} {Z}_s\geq 0)= \sqrt{\frac{2}{\sigma^2\pi}}$ and by (\ref{bpart}), for any $\epsilon>0$, there exists $K>0$ large enough such that
\begin{equation}
	\mathbb{E}_{a+g_l}\left(   \frac{| {Z}_t|}{\sqrt{t}}  \1_{\{   \frac{| {Z}_t|}{\sqrt{t}}\geq K  \}}  \big|  \inf_{s\leq t} {Z}_s \geq 0 \right)+ \E\big( m_1 \1_{\{ m_1\geq K\}}\big) \leq \epsilon,
\end{equation}
where we recall that $(m_s)_{s\in [0,1]}$ is a Brownian meander. Moreover by applying (\ref{apart}) to the function $G: (g_s)_{s\in [0,1]}\mapsto g_1 F((g_s)_{s\in [0,1]}) \1_{\{ g_1\leq K\}}$ we have
\begin{eqnarray}
	\nonumber \limsup_{t\to\infty} |  \mathbb{E}_{a+g_l}\left(  \frac{{Z}_t}{\sqrt{t}}  F\left\{  (t+l)^{-1/2}\Big(   g_{st}\1_{\{ s\leq \frac{l}{t} \}}   +{Z}_{st}    \1_{\{ s\geq \frac{l}{t} \}}\Big)_{s\in [0,1]} \right\}\big| \inf_{s\leq t} {Z}_s\geq 0 \right)
	- \mathbb E\big(m_1F((\sigma m_s)_{s\in [0,1]})\big)  |
	\\
	\leq   \limsup_{t\to\infty} \mathbb{E}_{a+g_l}\left(   \frac{| {Z}_t|}{\sqrt{t}}  \1_{\{   \frac{| {Z}_t|}{\sqrt{t}}\geq K  \}}  \big|  \inf_{s\leq t} {Z}_s \geq 0 \right)+ \mathbb E\big( m_1 \1_{\{m_1\geq K\}}\big) \leq \epsilon .
\end{eqnarray}
Recall that the under the probability $
\sqrt{\frac{2}{\pi}} m_1\cdot \mathbb{P}$, the process $(m_s)_{s\in [0,1]}$ has the same law as the three dimensional Bessel process $(R_s)_{s\geq 0}$. So by letting $\epsilon$ going to $0$, we get  
\begin{eqnarray*}
	&&\lim_{t\to\infty}  \mathbb{E}\left(   ({Z_t}+a+g_l )\1_{\{ \zeta^{a+g_l}>t  \}}\times\right.
	F\left\{ (t+l)^{-1/2}\Big( g_{st}\1_{\{ s\leq \frac{l}{t} \}}   +Z_{st}\1_{\{ s\geq \frac{l}{t} \}}\Big)_{s\in [0,1]} \right\} \Big)_{\big| (g_s)_{s\leq l}=(Z_i(s))_{s\leq l} }\\
		&&\qquad \qquad\qquad \qquad  \qquad\qquad  \qquad \qquad \qquad\qquad \qquad= (a +{Z_i(l)}) \mathbb E\Big(F((\sigma R_s)_{s\in [0,1]})\Big),
\end{eqnarray*}
and deduce that 
\begin{eqnarray*}
	&&\lim_{t\to\infty}\E\left(  \sum_{i=1}^\infty (a+Z_i(t+l))\ee^{-(1+\bar{p})Z_i(t+l)}\1_{\{{\zeta^a_i}> t+l \}} F\{(t+l)^{-1/2}(Z_i(u(t+l))_{u\in [0,1]} \} \big|\mathcal{F}_l\right)
	\\
	&&=\sum_{i=1} (a + Z_i(l))\ee^{-(1+\bar{p})Z_i(l)}\1_{\{{\zeta^a_i}>l \}} \mathbb E\Big(F((\sigma R_s)_{s\in [0,1]})\Big) \\
	&& =   \mathbb E\Big(F((\sigma R_s)_{s\in [0,1]})\Big) M'^{(a)}(l) .
\end{eqnarray*}
This concludes the proof of \ref{eqConditionalii}.
\hfill$\Box$

\subsection{$L^p$ bound and Concentration}

In this section we prove  Lemma \ref{Concentration}. {\it As there are a lot of estimates in this section, for notational convenience, the reader will note that $c$ always denotes a strictly positive constant which may vary in its value from line to line.} We start with a technical lemma. 
\begin{Lemma}
	\label{Momentdomination}Fix $a>0$.
	Under  (\ref{e.kappa}), (\ref{Asump2}) and (\ref{Asump1}). Let $\gamma\in (1,2)$ such that (\ref{integrab}) holds and $\epsilon \in (0,1)$. There exists $c>0$ such that for any $t\geq 1,\, a\leq t^\epsilon$,  we have 
	\begin{eqnarray}
		\label{eqMomentdomination}
		\E\left((\sqrt{t}W^{(a)}(t,\bar{p}))^\gamma \right) &\leq &c     \ee^{(\gamma-1)(1+\bar{p})a},
		\\
		\label{eqMomentdomination2}
		\E\left((M'^{(a)}(t))^\gamma\right) &\leq & c\ee^{(\gamma-1)(1+\bar{p})a}.
	\end{eqnarray}
	
\end{Lemma}
\noindent{\it Proof of Lemma \ref{Momentdomination}.} 
Appealing to the  the spinal decomposition described in Proposition \ref{BerRouau} (see in particular Remark \ref{laremarke}), we can write
\begin{eqnarray*}
	\E\left((\sqrt{t}W^{(a)}(t,\bar{p}))^\gamma \right) &=& t^{\frac{\gamma}{2}} \E\left( \sum_{i=1}^\infty \ee^{-(1+\bar{p})Z_i(t)}   \1_{\{{\zeta^a_i}>t \}}     (W^{(a)}(t,\bar{p}))^{\gamma -1}   \right)
	\\
	&=& t^{\frac{\gamma}{2}}\E^{(\bar{p})} \Big(  \E
	\big(  \sum_{i=1}^\infty \frac{\ee^{-(1+\bar{p})Z_i(t)}}{W(t,\bar{p})}1_{\{ {\zeta^a_i}>t  \}} W^{(a)}(t,\bar{p})^{\gamma-1}  \big| \mathcal{F}(t) \big)       \Big)
	\\
	&=& t^{\frac{\gamma}{2}}   \E^{(\bar{p})} \left( \1_{\{  \zeta_1^a>t \}} (W^{(a)}(t,\bar{p}))^{\gamma -1}   \right).
\end{eqnarray*}
Let $\mathcal{G}_1(\infty)$ be the sigma-field generated by the Poisson random measure  $N(dt, d\pi, k)$ to the fiber $\r_+\times \mathcal{P}\times \{1\}$. As $\gamma-1\in (0,1)$, by conditioning on $\mathcal{G}_1(\infty)$, Jensen's inequality gives us
\begin{eqnarray*}
	\E\left((\sqrt{t}W^{(a)}(t,\bar{p}))^\gamma \right) 
	&\leq&    t^{\frac{\gamma}{2}}\E^{(\bar{p})}\left(     \1_{\{ \zeta^a_1>t\}} \E^{(\bar{p})}\Big(W^{(a)}(t,\bar{p})\big| \mathcal{G}_1(\infty)\Big)^{\gamma -1}         \right).
\end{eqnarray*} 
Following in the style of reasoning found in e.g. \cite{HH}, to bound this conditional expectation we need to decompose it in a spine term and a additive term.  More precisely, recall that, if  $(r, \pi(r), 1)$ is a point in the Poisson point process $N(dt, d\pi, k)$, then the block $\Pi_1(r-)$ splits into $\pi(r)_{\big|\Pi_1(r-)}$, and the block after the split which contains $1$ is $\Pi_1(r)= \pi_1(r)\cap\Pi_1(r-)$. It follows that, there is  some index $j\geq 2$ such that $\Pi_i(t)\subset \pi_j(r)\cap \Pi_1(r-)$, where $\pi_j(r)$ stands for the $j-$th block of the partition $\pi(r)$. In other words, we can write
\begin{eqnarray}
	\nonumber    && \E^{(\bar{p})}\Big(W^{(a)}(t,\bar{p})\big| \mathcal{G}_1(\infty)\Big)= (\mathtt{S}) + (\mathtt{A}),
\end{eqnarray}
where $(\mathtt{A})$ and $(\mathtt{S})$ are respectively called the ``additive term" and the ``spine term" and defined by: 
\begin{eqnarray}
	\label{SmoinsA} &&(\mathtt{A}) := \sum_{r\in \mathcal{D}_1\cap [0,t]}\ee^{ \Phi(\bar{p})r } \sum_{j=2}^\infty  |\pi_j(r)\cap \Pi_1(r-)|^{1+\bar{p}}  
	\P^{(\bar{p})}\left( \zeta^{a+x}>t-r \right)_{\big| x= -\log|\pi_j(r)\cap \Pi_1(r-)|   - r\Phi'(\bar{p}) } ,
	\\
	&&(\mathtt{S}):= \ee^{-(1+\bar{p})Z_1(t)}. \label{SplusA} 
\end{eqnarray}
We first bound the ``spine term" $(\mathtt{S})$. Recall that $\Phi(\bar{p})=(1+\bar{p})\Phi'(\bar{p})$, then  we can write
\begin{eqnarray*}
	\E^{(\bar{p})}\left( \1_{\{  \zeta^a_1>t \}} (\mathtt{S})^{\gamma-1}  \right) 
	&=& \mathbb{E}\left(    \1_{\{ \zeta^a>t \}} \ee^{ -(\gamma-1)(1+\bar{p}){Z}_t}     \right),
\end{eqnarray*}
By localizing the value of $Z_t$, then using (\ref{2.6}) we obtain
\begin{eqnarray}
	\nonumber 	\mathbb{E}\left(    \1_{\{  \zeta^a_1>t \}} \ee^{ -(\gamma-1)(1+\bar{p}){Z}_1(t)}     \right) &\leq &  \ee^{(1+\bar{p})(\gamma-1) a}  \sum_{k\geq 0} \ee^{ -(\gamma-1)(1+\bar{p}) k} \mathbb{P}_a\left(    \inf_{s\leq t} {Z}_s \geq 0,\,   {Z}_t\in [k,k+1]   \right)
	\\
	\label{Sbound}	&\leq &\ee^{(1+\bar{p})(\gamma-1) a}  \frac{(1+a)}{t^{\frac{3}{2}}} \notag\\
	&\leq& c \ee^{(\gamma-1)(1+\bar{p})a},
\end{eqnarray}
where we used $a\leq t^{\epsilon}$ in the last inequality. To prove inequality (\ref{eqMomentdomination}) it remains to show that there exists $c>0$ such that
\begin{equation}
	\label{equivproof} t^\frac{\gamma}{2} \E^{(\bar{p})}\left( \1_{\{\zeta^a_1>t  \}}  (\mathtt{A})^{\gamma-1} \right)\leq c \ee^{(\gamma- 1)(\bar{p}+1)a}.
\end{equation}
By using \eqref{2.5}, for any $x\geq -a$  there exists a constant  $c>0$ such that  
\begin{equation}
	\P^{(\bar{p})}\left(\zeta^{x+a} >t-r \right)\leq c \frac{ (1+x+a) }{ (1+t-r)^{\frac{1}{2}}}.
\end{equation}
Plugging this inequality into (\ref{SmoinsA}) and using the identity $|\pi_j(r)\cap \Pi_1(r-)|=|\pi_j(r)||\Pi_1(r-)| $ yields
\begin{eqnarray*}
	(\mathtt{A})&\leq& c   \sum_{r\in \mathcal{D}_1\cap [0,t]} \frac{  \ee^{ -(1+\bar{p})Z_{1}(r-) }}{(1+t-r )^{\frac{1}{2}}}    \sum_{j=2}^\infty  |\pi_j(r)|^{1+\bar{p}}  \left[1 -\log|\pi_j(r)| + Z_{1}(r-)   +a \right].
\end{eqnarray*}
It is easy to show that 
\begin{eqnarray*}
	\sum_{j=2}^\infty  |\pi_j(r)|^{1+\bar{p}}  \Big[1 -\log|\pi_j(r)|  + Z_{1}(r-)+  a \Big]  \leq   (  Z_{1}(r-)+a)
	\left( \sum_{j=2}^\infty \big[1-\log  |\pi_j(r)|   \big] | \pi_j(r)|^{ 1+\bar{p}}\right).
\end{eqnarray*} 
Moreover for any $x,y\geq 0$,  $(x +y)^{\gamma-1}\leq x^{\gamma-1} +y^{\gamma-1}$, recall also that $\Phi(\bar{p})=(1+\bar{p})\Phi'(\bar{p})$,  it yields
\begin{eqnarray*}
	(\mathtt{A})^{\gamma-1}   \leq    c \sum_{r\in \mathcal{D}_1\cap [0,t]} \frac{\ee^{-(\gamma-1)( 1+\bar{p})Z_{1}(r-)}}{(1+t-r)^{\frac{(\gamma-1)}{2}}}  (  Z_{1}(r-)+a)^{(\gamma-1)}
	\left( \sum_{j=2}^\infty [1-\log  |\pi_j(r)|   ] | \pi_j(r)|^{ 1+\bar{p}}\right)^{(\gamma-1)}  .
\end{eqnarray*}

 Recalling the conclusion of Proposition \ref{BerRouau} (ii) and again appealing to \eqref{2.5}, we have that, for any $r\in [0,t]$, there exists a constant $c$ such that 
\begin{eqnarray*}
\P^{(\bar{p})}\left( \zeta^a_1>t  | \mathcal{F}(r) \right)&& \leq 
c\frac{(1+a+ Z_1(r-) -\log |\pi_1(r)| )  )}{(1+t-r)^{1/2}} .
\end{eqnarray*}
So by computing the predictable compensator of such additive functional, it leads to
\begin{eqnarray}
	\nonumber 	\E\left( \1_{\{ \zeta^a_1>t  \}}  (\mathtt{A})^{\gamma-1} \right)
	&\leq&  \mathbb{E}\left(     \int_0^t \1_{\{   \zeta^a>s  \}}\frac{( a + Z_{s-})^{\gamma } \ee^{-(\gamma-1)( 1+\bar{p})Z_{s-}}     }{(1+t-s)^{\frac{\gamma}{2}}}        ds \right) 
	\\	
	\label{comoing}	&& \qquad\qquad\qquad\qquad\qquad\qquad \times  \int_{\mathcal{P}}\mu(d\pi )|\log |\pi_1|| |\pi_1|^{\bar{p}}    \left( \sum_{j=2}^\infty [1-\log  |\pi_j|   ] | \pi_j|^{ 1+\bar{p}}\right)^{(\gamma-1)}.
\end{eqnarray}
The crucial point here is that 
\begin{equation}
	\label{Crucia}
	\int_{\mathcal{P}}\mu(d\pi )  |\log |\pi_1|| |\pi_1|^{\bar{p}}    \left( \sum_{j=2}^\infty [1-\log  |\pi_j|   ] | \pi_j|^{ 1+\bar{p}}\right)^{(\gamma-1)}<\infty.
\end{equation}
Indeed by the Proposition (\ref{BerRouau}) and the identity (3) of \cite{HKK}, we have that
\begin{eqnarray*}
	&&\int_{\mathcal{P}}\mu(d\pi ) |\log |\pi_1|| |\pi_1|^{\bar{p}}    \left( \sum_{j=2}^\infty [1-\log  |\pi_j|   ] | \pi_j|^{ 1+\bar{p}}\right)^{(\gamma-1)}
	\\
	&&\qquad\qquad \qquad =\int_{\mathcal{S}}\nu(d{\bf s} )\left(  \sum_{i=1}^\infty (-\log s_i) s_i^{1+\bar{p}}\right)   \left( \sum_{i=2}^\infty [1-\log  s_i   ] s_i^{ 1+\bar{p}}.\right)^{(\gamma-1)}
\end{eqnarray*}
Trivially $ \sum_{i=1}^\infty s_i^{1+\bar{p}}\leq  \sum_{i=2}^\infty s_i^{1+\bar{p}} +1  $. Moreover thanks to the inequality $(x+y)^{\gamma-1}\leq x^{\gamma-1} +y^{\gamma- 1}$ (recall that $\gamma-1\in (0,1)$),  we also have, for some   small $\epsilon>0$ and constant $c = c(\epsilon)>0$,
\begin{eqnarray*}
	\left(  \sum_{i=1}^\infty (-\log s_i)s_i^{1+\bar{p}}\right)   \left( \sum_{i=2}^\infty [1-\log  s_i   ] s_i^{ 1+\bar{p}}\right)^{(\gamma-1)}	&\leq &   c \left\{  \left(  \sum_{i=2}^\infty s_i^{1+\bar{p}-\epsilon}\right)^{\gamma} + \left(  \sum_{i=2}^\infty s_i^{1+\bar{p}-\epsilon}\right)^{\gamma-1}  \right\}   .
\end{eqnarray*}
 The term on the right-hand side above is integrable with respect to $\nu(d{\bf s})$ by (\ref{holder2}), from which (\ref{Crucia}) now follows. Coming back to (\ref{comoing}), we are now in position to conclude the proof of (\ref{equivproof}). First we apply the Fubini-Tonelli Theorem and the inequality (\ref{Crucia}), which yields
\begin{eqnarray*}
  t^\frac{\gamma}{2} \E\left( \1_{\{\zeta_1^a>t\}}  (\mathtt{A})^{\gamma-1} \right)
&\leq &c  t^{\frac{\gamma}{2}}  \int_0^t  \mathbb{E}\left(\1_{\{ \zeta^a>s\}}   \frac{(  a + Z_{s-})^{\gamma } \ee^{-(\gamma-1)( 1+\bar{p})Z_{s-}}     }{(t-s+1)^{\frac{\gamma}{2}}}    \right)    ds.
	\\
	&\leq & c \ee^{(\gamma-1)( 1+\bar{p})a}  \left\{  \int_{0}^{\frac{t}{2}}   \mathbb{E}_a\left( \1_{\{ \zeta^0>s \}}  \ee^{-\frac{(\gamma-1)}{2}( 1+\bar{p}){Z}_s}      \right) ds \right.\\
	&&\left.\hspace{3cm}+   t^{\frac{\gamma}{2}} \int_{\frac{t}{2}}^t \frac{       \mathbb{E}_a\left( \1_{\{ \zeta^0>s \}}   \ee^{-\frac{(\gamma-1)}{2}( 1+\bar{p}){Z}_s} \right)      }{(1+t-s)^{\frac{\gamma}{2}} }ds  \right\},
\end{eqnarray*}
where we used in the last line the inequality $x^{\gamma}\ee^{-a x} \leq c(a) \ee^{-\frac{a}{2}x}$, $\forall a\geq 0,\, x\geq 0$. By  Lemma \ref{PourRW2}, the first term is smaller than $\ee^{(\gamma-1)(1+\bar{p})}c(\frac{\gamma-1}{2}(1+\bar{p}))$. By Lemma \ref{2ineq} and a localization of $Z_t$, for any $s\in [\frac{t}{2},t]$ we have
$$ \Pp_a\left( \1_{\{  \zeta^0>s \}}   \ee^{-\frac{(\gamma-1)}{2}( 1+\bar{p}){Z}_s} \right)   \leq (1+a)t^{-\frac{3}{2}}.$$ 
It therefore follows  that 
$$  \int_{\frac{t}{2}}^t t^{\gamma\over 2} \mathbb{E}_a\left( \1_{\{ \zeta^0>s \}}   \ee^{-\frac{(\gamma-1)}{2}( 1+\bar{p}){Z}_s} \right)  ds     \leq   t^{\gamma\over 2} t^{{-\gamma\over 2} +1} (1+a)t^{-\frac{3}{2}}\leq  \frac{(1+a)}{\sqrt{t}}  \leq 1, $$
where we have used $a\leq t^\epsilon$ in the last inequality. We deduce that
\begin{equation*}
	t^\frac{\gamma}{2} \E\left( \1_{\{\zeta_1^a>t  \}}  (\mathtt{A})^{\gamma-1} \right)\leq c \ee^{(\gamma-1)( 1+\bar{p})a},
\end{equation*}
which ends the proof of the inequality (\ref{eqMomentdomination}). The proof (\ref{eqMomentdomination2}) is very similar. Indeed it suffices to follow the same scheme, we just point out here the main arguments. 
We start by using  Proposition \ref{BerRouau} and write
$$ \E\left(  (M'^{(a)}(t))^\gamma  \right)= \E^{(\bar{p})}\left(  \1_{\{  \zeta^a_1>t \}}(a  + Z_1(t)) (M'^{(a)}(t))^{\gamma-1}  \right) $$
We then take  conditional expectation with respect to $\mathcal{G}_1(\infty)$ then use Jensen's inequality. Similarly we can decompose
\begin{equation}
	\E^{(\bar{p})}\left(  M'^{(a)}(t)   \big| \mathcal{G}_1(\infty) \right)= (\mathtt{S})+(\mathtt{A}),
\end{equation}
with
\begin{eqnarray*}
	&&(\mathtt{A}) := \sum_{r\in \mathcal{D}_1\cap [0,t]}\ee^{ \Phi(\bar{p})r } \sum_{j=2}^\infty  |\pi_j(r)\cap \Pi_1(r-)|^{1+\bar{p}}  (\xi_{r-}-r^-\Phi'(\bar{p})+a) 
	\\
	\nonumber &&\qquad\qquad\qquad\qquad  \times\mathbb{E}\left(  (a+ Z_{t-r})  \1_{\{\zeta^{x+a}>t-r\}} \right)_{\big| x= -\log|\pi_j(r)\cap \Pi_1(r-)|   - r\Phi'(\bar{p}) } ,
	\\
	&&(\mathtt{S}):= \ee^{\Phi(\bar{p})t - (1+\bar{p})\xi_t} 
\end{eqnarray*}
With regard to the ``spine term", by using (\ref{2.6}) we have
\begin{equation*}
	\E^{(\bar{p})}\left( \1_{\{ \zeta_1^a>t \}} (a+Z_1(t)) (\mathtt{S})^{\gamma-1}  \right) \leq   \frac{c(1+a)   \ee^{(\gamma-1)(1+\bar{p}) a} }{t^{\frac{3}{2}}}\leq c' \ee^{(\gamma-1)(1+\bar{p}) a}
\end{equation*}
For the ``additive term",  for any $j\geq 2$, by (\ref{toadd}), the expectations terms are bounded by $(1-\log |\pi_j(r)|+ Z_1(r-))$. Then by some elementary inequalities we can get
\begin{eqnarray*}
	(\mathtt{A})^{\gamma-1} \leq   c \sum_{r\in \mathcal{D}_1\cap [0,t]} \frac{\ee^{-( 1+\bar{p})(\gamma-1)Z_1(r-)}}{(1+t-r)^{\frac{(\gamma-1)}{2}}} 
(  a + Z_1(r-))^{2(\gamma-1)}
	\left( \sum_{j=2}^\infty [1-\log  |\pi_j(r)|   ] | \pi_j(r)|^{ 1+\bar{p}}\right)^{(\gamma-1)}  .
\end{eqnarray*}
Finally  by using the same type of arguments (Markov property with (\ref{toadd}), computation of the predictable compensator, Fubini-Tonelli, ) we can similarly obtain
\begin{eqnarray*}
	\E\left( \1_{\{\zeta_1^a>t  \}}  (a+{Z_1(t)}) (\mathtt{A})^{\gamma-1} \right) \leq    c\ee^{(\gamma-1)( 1+\bar{p})a}  \int_0^t   \mathbb{E}_a\left( \1_{\{  \inf_{u\leq s} {Z}_s \geq 0 \}}  {Z}_s ^{ \gamma  } \ee^{-(\gamma-1)( 1+\bar{p}){Z}_s}      \right) ds  
	\\
	\times \int_{\mathcal{P}}\mu(d\pi )|\log |\pi_1| | |\pi_1|^{\bar{p}}    \left( \sum_{j=2}^\infty [1-\log  |\pi_j|   ] | \pi_j|^{ 1+\bar{p}}\right)^{(\gamma-1)},
\end{eqnarray*}
where the first and the second terms are finite respectively thanks to (\ref{PourRW2}) and  (\ref{Crucia}). This concludes the proof of Lemma \ref{Momentdomination}.\hfill$\Box$\\

Now we are in position to prove the Proposition \ref{Concentration}. This part appeals to classical $L^p$ estimates;  see for instance \cite{HH} or \cite{EHK10}.

\bigskip

\noindent{\it Proof of Proposition \ref{Concentration}.} We recall a powerful inequality proved in Lemma 1 of \cite{Big92}. For independent, zero mean random variables $\{Z_1,..., Z_n\}$ and $p\in [1,2]$, we have
\begin{equation} 
	\label{russkov} \E\left( \left|\sum_{i=1}^n Z_i\right|^p\right) \leq 2^p\sum_{i=1}^n \E(|Z_i|^p).
\end{equation}
We shall use this inequality with $\gamma\in (1,2)$ such that (\ref{integrab}) holds. An application of  the Markov inequality and then  the fragmentation property gives us
\begin{eqnarray*}
	\lefteqn{\P\left(  \left| \E\left( (t+l)^{1/2}\, W^{(a)}(t+l,\bar{p}) \big|\mathcal{F}_l\right)-   (t+l)^{1/2}\, W^{(a)}(t+l,\bar{p})\right|\geq \epsilon \right)}&&\\
&&	\leq
	\epsilon^{-\gamma} \E\left( \left| \sum_{i=1}^\infty \ee^{-(1+\bar{p})Z_i(l)} \1_{\{{\zeta^a_i}>l \}}   (t+l)^{1/2}\left[  \E\left( W^{(a+x)}(t,\bar{p})  \right) - W^{(a+x)}(t,\bar{p})\right]_{\big| x={Z_i(l)}} \right|^\gamma \right)
\end{eqnarray*}
Applying the inequality (\ref{russkov}), we deduce that there exists $c>0$ such that
\begin{eqnarray*}
	\lefteqn{\P\left(  \left| \E\left( (t+l)^{1/2}\, W^{(a)}(t+l,\bar{p}) \big|\mathcal{F}_l\right)-   (t+l)^{1/2}\, W^{(a)}(t+l,\bar{p})\right|\geq \epsilon \right)}&&
	\\
	&& \leq c\epsilon^{-\gamma} \E\left(  \sum_{i=1}^{\infty} \ee^{-(1+\bar{p})\gamma Z_i(l)}  \1_{\{{\zeta^a_i}>l \}}\E\left(\left| (t+l)^{1/2} W^{(a+x)}(\bar{p},t)\right|^\gamma \right)_{\big| x={Z_i(l)}}\right) 
	\\
	&&\leq c  \frac{ (t+l)^\frac{\gamma}{2}      
	}{t^{\frac{\gamma}{2}}\epsilon^\gamma  } \E\left( \sum_{i=1}^\infty \ee^{-(1+\bar{p})\gamma Z_i(l)} \1_{\{{\zeta^a_i}>l \}} \ee^{(\gamma-1)(1+\bar{p})(a+ {Z_i(l)}) }  \right),\\
	&&=  c\frac{ (t+l)^\frac{\gamma}{2}    \ee^{(\gamma-1)(1+\bar{p}) a} }{t^{\frac{\gamma}{2}}\epsilon^\gamma  }   \E\left(   \sum_{i=1}^\infty  \ee^{-(1+\bar{p})Z_i(l)}  \1_{\{{\zeta^a_i}>l \}}   \right) .
\end{eqnarray*}
where in the penultimate inequality we have used (\ref{eqMomentdomination}).
It remains to use  Corollary \ref{Manytoone} then the inequality in (\ref{2.5}) to conclude that
\begin{eqnarray*}
	\lefteqn{\hspace{-3cm}\P\left(  \left| \E\left( (t+l)^{1/2}\, W^{(a)}(t+l,\bar{p}) \big|\mathcal{F}_l\right)-   (t+l)^{1/2}\, W^{(a)}(t+l,\bar{p}) \right|\geq \epsilon \right)
	} &&\\
	&&\leq   c  \frac{ (t+l)^\frac{\gamma}{2}   \ee^{(\gamma-1)(1+\bar{p}) a} }{t^{\frac{\gamma}{2}}\epsilon^\gamma  } \mathbb{P}\left(  \zeta^a>l  \right)
	\\
	&&  \leq\frac{c(t+l)^\frac{\gamma}{2}}{t^{\frac{\gamma}{2}}\epsilon^\gamma  } \frac{ (1+a)    \ee^{(\gamma-1)(1+\bar{p}) a} }{\sqrt{l}},
\end{eqnarray*}
which ends the proof of  (\ref{SeneCond}). The proof of (\ref{BesCond}) is very similar and left to the reader. \hfill$\Box$

\section{BRW, BBM and log-correlated Gaussian field} 
\label{BRWBBMMC}
\subsection{Branching random walk}

The model of the branching random walk (BRW) can be described as follows. Initially,  a single particle denoted $\emptyset$ sits at the origin.  Its children together with their  displacements,  form a point process $\Theta$ on $\r$ and  the first generation of the branching random walk.  These children have children of their own which form the second generation, and behave, relative to their respective positions at birth, like independent copies of the same point process $\Theta$, and so on.

Let $\mathbb{T}$ be the genealogical tree of the particles in the branching random walk. Plainly, $\mathbb{T}$ is a Galton-Watson tree.  We write $|z|=n$ if a particle  $z$  is in the $n$-th generation, and denote its position by $V(z)$ ($V(\emptyset)=0$). The collection of positions $(V(z),z\in \mathbb{T})$ is our branching random walk.

We assume throughout the paper the following conditions
\begin{eqnarray}
	\label{chap3criticalcondition1}
	\E\Big(\underset{|x|=1}{\sum}  \ee^{-V(x)} \Big)=1,&   &\E\Big(   \underset{|x|=1}{\sum}1  \Big) >1 ,\quad \text{  and }
	\\
	\label{chap3criticalcondition2} \E\Big(\underset{|x|=1}{\sum}V(x)\ee^{-V(x)} \Big)=0,&   \qquad &\sigma^2:=\E \Big(\underset{|x|=1}{\sum}V(x)^2\ee^{-V(x)} \Big)<\infty.
\end{eqnarray}
The branching random walk is then said to be in the boundary case (Biggins and Kyprianou \cite{Bky05}). Moreover we will assume that for some $\gamma\in (0,1)$,
\begin{eqnarray}
	\label{Pmoment}
	\E\Big( \sum_{|z|=1}\ee^{-V(u)} \Big(   \sum_{|z|=1}  V(u)_+\ee^{-V(u)} \Big)^\gamma \Big)<\infty,
\end{eqnarray}
where, for any $x\in \r$, $x_+= \max(0,x)$. Let us introduce
\begin{equation*}
	W_{n}:=\underset{|x|=n}{\sum}\ee^{-  V(x)},\qquad n\geq 1,
\end{equation*} 
the {\it critical additive martingale} associated with the branching random walk and 
$$M_n:=\underset{|u|=n}{\sum}V(u)\ee^{-V(u)}, \qquad n\ge 1.$$
its {\it critical derivative martingale}. On the set of non-extinction ($\mathbb{T}$ is infinite), by Biggins and Kyprianou \cite{BKy04} we know that $M_n$ admits an almost sure limit which we will denote by $M_\infty$. The main result of this section is to give a different proof of the following result, which first appeared in \cite{AShi11}.
\begin{theorem}
	Assume (\ref{chap3criticalcondition1}), (\ref{chap3criticalcondition2}) and (\ref{Pmoment})  we have
	\begin{equation}
		\label{BRWtheorem}
		\underset{n\to\infty}{\lim} \sqrt{n} \frac{W_n}{M_n}=\left(\frac{2}{\pi\sigma^2}\right)^{\frac{1}{2}},\qquad \text{ in } \P^*\text{ probability}.
	\end{equation}
	where 	$\P^*(\cdot):=\P\left(\cdot\, |\, \mathbb{T} \text{ is infinite}\right)$. 
\end{theorem}
Note, by several simple  convexity inequalities, it can be easily checked that (\ref{Pmoment}) is stronger than the assumptions used for the corresponding version of the above theorem  in \cite{AShi11}. Thus the emphasis of this section this lies more in the direction of the method of proof, which is relatively shorter. 

In the following section we  briefly collect some preliminary results on the branching random walk (change of probabilities, an associated one-dimensional random walk),  which are lifted in their present form from    \cite{AShi11}. 

\subsection{The many-to-one Lemma}
Let $(V(x))_{x\in\mathbb{T}}$ be a branching random walk starting at the origin and satisfying (\ref{chap3criticalcondition1}) and (\ref{chap3criticalcondition2}). Let $(S_n)_{n\geq 0}$ be the random walk such that $S_0=0$ and the law of $(S_1)$ is given by 
\begin{equation}
	\label{chap3defSn}
	\E\left(f(S_1)\right):=\E \Big(\underset{|z|=1}{\sum}f(V(z))\ee^{-V(z)} \Big),
\end{equation}
for all positive, bounded and measurable functions.
The conditions (\ref{chap3criticalcondition1}) and (\ref{chap3criticalcondition2}) imply that $(S_n)_{n\geq 0}$ is a mean zero random walk and $\E(S_1^2)=\sigma^2<\infty$. A straightforward inductive argument shows that, for any $n\geq 0$ and $g:\r^n\to \r^+$ bounded and measurable, we have{\string:}
\begin{equation}
	\label{chap32.1}
	\E \Big(\underset{|x|=n}{\sum}g\left(V(x_1),...,V(x_n)\right) \Big)=\E\left(\ee^{S_n}g(S_1,...,S_n)\right),
\end{equation}
where for any $i=1,\cdots, n$, $x_i$ is the unique ancestor of $x$ such that $|x_i|=i$. Equality (\ref{chap32.1}) is the so-called {\it many-to-one identity} which plays a fundamental role in many computations relating to the linear semi-group of the branching random walk. The presence of the random walk $(S_n)_{n\geq 0}$ was first explained in  Lyons \cite{Lyo97} and Biggins, Kyprianou \cite{BKy04}. Let $(\mathcal{F}_n)_{n\geq 0}$ be the sigma-algebra generated by the branching random walk in the first $n$ generations. We call a {\it spine} in $\mathbb{T}$ an infinite line of descent starting from the root. Since the landmark work of Lyons, Pemantle and Peres \cite{LPP95} and Lyons \cite{Lyo97}, the {\it spine decomposition} is a widespread technique to study the branching random walk. To describe it, we need to introduce the change of measure
$${\Q}_{\big|\mathcal{F}_n}:= W_{n} . \P_{\big|\mathcal{F}_n}, \qquad n\geq 0.$$ 
We denote by $\hat{\Theta}$ the point process with Radon-Nikodym derivative $\sum_{|x|=1}\ee^{- V(x)}$ with respect to the law of $\Theta$.  In \cite{Lyo97} we can find the following description of branching random walk under $\Q${\string:}

\begin{itemize}
	\item[(i)] $w_0=\emptyset$ gives birth to particles distributed according to $\hat{\Theta}$.
	\item[(ii)] Choose $w_1$ among children of $w_0$ with probability proportional to $\ee^{- V(x)}$.
	\item[(iii)] $\forall n\geq 1$, $w_n$ gives birth to particles distributed according to $\hat{\Theta}$ (with $u=V(w_n)$).
	\item[(iv)] Choose $w_{n+1}$ among the children of $w_n$ with probability proportional to $\ee^{-V(x)}$.
	\item[(v)] Subtrees rooted at all other brother particles are independent branching random walks with the same distribution as under $\P$.
\end{itemize}
\paragraph
\noindent From this description we deduce the three followings facts:{\string:}
First,  $\Q(\text{non-extinction})=1$. Second, 
for any $n$ and any vertex $x$ with $|x|=n$, we have
	\begin{equation}
		\label{Spiine}	\Q(w_n=x|\mathcal{F}_n)= \frac{\ee^{- V(x)}}{W_n}.
	\end{equation}
Finally, the spine process $(V(w_n)), n\geq 0)$ under $\Q$, is distributed as a standard cantered random walk.

\subsection{The renewal function and random walk asymptotics}
The random walk $(S_n)_{n\geq 0}$ is centred with  $\sigma^2=\E[S_1^2]\in(0,\infty)$. Let $h_0$ be its renewal function defined by
\begin{equation}
	\label{definih0}
	h_0(u):=\underset{j\geq 0}{\sum}\P\left( \underset{i\leq j-1}{\min}S_i>S_j\geq -u  \right),\qquad u\geq 0.
\end{equation}
Note that $h_0$ is increasing and $h_0(0)=1$. For any $u\geq 0$, $h_0$ satisfies
\begin{equation}
	\label{chap33.40}
	h_0(u)=\E\left(h_0(S_1+u)\1_{\{S_1\geq -u\}}\right).
\end{equation}
Write  $ \underline{S}_n:=\underset{j\leq n}{\min}\,\, S_j,\, \forall   n\geq 0$, for the running minimum. It is known that there exists $c_0>0$ and $\theta>0$ such that
\begin{equation}
	\label{chap3lldda}
	c_0:=\underset{u\to\infty}{\lim}\frac{h_0(u)}{u},\qquad \P\left(\underline{S}_n\geq -u\right)\underset{n\to\infty}{\sim}\frac{1}{c_0}\sqrt{\frac{2}{\pi \sigma^2}}\frac{ h_0(u)}{n^\frac{1}{2}},\qquad \forall u\geq 0.
\end{equation}
As in \cite{AShi11}, we will need the following uniform version of (\ref{chap3lldda}){\string:}  as $n\to \infty$,
\begin{equation}
	\label{chap3lldda2}
	\P\left(\underline{S}_n\geq -u\right)=  \frac{1}{c_0} \sqrt{\frac{2}{\pi \sigma^2}}\frac{ h_0(u)+ o(1)}{n^\frac{1}{2}},
\end{equation}
uniformly for $u\in [0, (\log n)^{30}]$. Finally we mention an inequality proved in \cite{AShi11}. For $u>0,a\geq0,\,b\geq 0$ and $n\geq 1$,
 there exists $c>0$ such that 
 	\begin{equation}
		\label{chap3Aine}
		\P\left(\underline{S}_n\geq -a,\,b-a\leq S_n\leq b-a+u\right)\leq c\frac{(u+1)(a+1)(b+u+1)}{n^{\frac{3}{2}}}.
	\end{equation}

	\subsection{Proof of Theorem \ref{BRWtheorem}}
Once again, on account of the number of estimates that are made in this section, the reader notes that $c$ is reserved for a generic positive constant that my vary from line to line.
	
	For any $u\in \mathbb{T}$, let $\underline{V}(u):= \min_{i=1,\cdots, n}V(u_i)$, where $u_i$ is the unique ancestor of $u$ such that $|u_i|=i$. We introduce the truncated {\it additive} and {\it derivative} martingales, 
	\begin{equation}
		W_n^{(a)}:= \sum_{|z|=n} \ee^{-V(z)} \1_{\{ \underline{V}(z)\geq -a  \}},\qquad M_n^{(a)}:= \sum_{|z|=n} h_0(a+V(z)) \ee^{-V(z)} \1_{\{ \underline{V}(z)\geq -a  \}}
	\end{equation}
	The proof of the Theorem \ref{BRWtheorem} follows the same scheme as the proof of Theorem \ref{mainresult}. Proposition \ref{BRWlargenumber} and Proposition \ref{Concentration} below are the two main steps. 
	\begin{Proposition}
		\label{BRWlargenumber}
		For any $a,\, l>0$, we have
		\begin{equation}
			\lim_{n\to\infty} \sqrt{n} \E( W_{n+l}^{(a)}\big| \mathcal{F}_l)\overset{\text{a.s}}{=} \frac{1}{c_0}{\sqrt{\frac{2}{\pi \sigma^2}}M^{(a)}_l}.
		\end{equation}
	\end{Proposition}
	\noindent{\it Proof of Proposition \ref{BRWlargenumber}.} By the branching property we have
	\begin{eqnarray*}
		\E\left( \sqrt{n} W_{n+l}^{(a)}\big|\mathcal{F}_l \right)&=& \sum_{|z|=l} \ee^{-V(z)} \1_{\{  \underline{V}(z)\geq -a\}} \E\left(  \sum_{|u|=n} \ee^{-V(u)} \1_{\{ \underline{V}(u)\geq -a+ V(z)  \}}   \right)
		\\
		&=&  \sum_{|z|=l} \ee^{-V(z)} \1_{\{  \underline{V}(z)\geq -a\}} \sqrt{n} \P\left( \underline{S}_n \geq -V(z) +a \right)
		\\
		&\underset{n\to \infty}{\longrightarrow} &\sum_{|z|=l} \ee^{-V(z)} \1_{\{  \underline{V}(z)\geq -a\}}\frac{1}{c_0} \sqrt{\frac{2}{\pi \sigma^2}} h_0(V(z)+a) = \frac{1}{c_0}\sqrt{\frac{2}{\pi \sigma^2}} M^{(a)}_l.
	\end{eqnarray*}
	where we have used the many-to-one Lemma then (\ref{chap3lldda}).  \hfill$\Box$

	\begin{Proposition}
		\label{ConcentrationBRW}
		Assume (\ref{chap3criticalcondition1}) and (\ref{chap3criticalcondition2}) and (\ref{Pmoment}). There exists $c>0$, $\gamma\in (1,2)$ such that:
		
		(i) For any $\epsilon >0,\, l,\, t\geq 0$ we have
		\begin{equation}
			\label{SeneCondBRW}
			\P\left(  \Big| \E\left( \sqrt{n+l}\, W_{n+l}^{(a)} \big|\mathcal{F}_l\right)-   \sqrt{n+l}\, W_{n+l}^{(a)}\Big|\geq \epsilon \right)\leq \frac{c(n+l)^\frac{\gamma}{2}}{t^{\frac{\gamma}{2}}\epsilon^\gamma \sqrt{l}}(1+a)\ee^{(\gamma-1) a}  .
		\end{equation}
	\end{Proposition}

	We start by proving the following Lemma:
	\begin{Lemma}
		\label{BRWPmoment}
		Under (\ref{chap3criticalcondition1}) and (\ref{chap3criticalcondition2}). Let $\gamma\in (1,2)$ such that (\ref{Pmoment}) holds. There exists $c>0$ such that for any $t\geq 1,\, a\leq t^{\epsilon}$, we have
		\begin{equation}
			\E\left(  \left( \sqrt{n} W^{(a)}_n \right)^\gamma \right)\leq c\ee^{(\gamma-1)a}.
		\end{equation}
	\end{Lemma}
	\noindent{\it Proof of Lemma \ref{BRWPmoment}.} By the representation given on the measure $\Q$,
	\begin{eqnarray}
		\E\left( (\sqrt{n} W_n^{(a)})^{\gamma} \right)=  n^{\frac{\gamma}{2}} \E_{\Q} \left( \1_{\{ \underline{V}(z)\geq -a \}} (W_n^{(a)})^{\gamma-1} \right)
	\end{eqnarray}
	Let us decompose $W_n^{(a)}$ using the spine. Specifically, we have 
	\begin{eqnarray*}
		\E\left( (\sqrt{n} W_n^{(a)})^{\gamma} \right)&= & n^{\frac{\gamma}{2}} \E_{\Q} \left( \1_{\{ \underline{V}(w_n)\geq -a \}} \Big( \ee^{-V(w_n)} + \sum_{k=1}^{n}\ee^{-V(w_{k-1})} \sum_{u\in \Omega(w_k)} \ee^{-[V(u)-V(w_{k-1})]} W_n^{(a)}(u)  \Big)^{\gamma-1} \right)
	\end{eqnarray*}
	with $\Omega(w_k):= \left\{|x|=k: w_{k-1}^{(\alpha)}\text{ is the parent of } x,\, x\neq w_k^{(\alpha)}\right\} $ and,  for all  $u\in \Omega(w_k)$, \[
	W_n^{(a)}(u):= \sum_{|z|=n,\, z\geq u} \ee^{-[V(z)-V(u)]} \1_{\{ \underline{V}(z)\geq -a \}}.\]
	 This decomposition leads us to introduce
	\begin{eqnarray}
		(\mathtt{S})&:=& \1_{\{ \underline{V}(w_n)\geq -a \}}  \ee^{- V(w_n)},
		\\
		(\mathtt{A})&=& \1_{\{ \underline{V}(w_n)\geq -a \}}  \sum_{k=1}^{n}\ee^{- V(w_{k-1})}   \sum_{u\in \Omega(w_k)} \ee^{-[V(u)-V(w_{k-1})]} W_n^{(a)}(u)   
	\end{eqnarray}
	For the ``spine term" $(\mathtt{S})$ by (\ref{chap3Aine}) we get
	\begin{equation}
		n^{\frac{\gamma}{2}}\E_\Q\left( (\mathtt{S})^{\gamma-1} \right)\leq \sum_{k\geq 0} n^{\gamma \over 2} \ee^{-(\gamma-1) k} \P_a( \underline{S}_n\geq 0,\, S_n\in [k,k+1]   )\leq \frac{\ee^{(\gamma-1) a} a}{n^{\frac{3-\gamma}{2}}}
	\end{equation}
	For $(\mathtt{A})$, the ``additive term", we introduce $
	\mathcal{G}_\infty:= \sigma(V(w_k),u, V(u),\, u\in \Omega(w_k),\, k\geq 0  )$ the sigma-field generated by the spine and its brothers. By conditioning according to $\mathcal{G}_\infty$ and using the Jensen inequality we get
	\begin{eqnarray*}
		\lefteqn{\E_\Q\left( (\mathtt{A})^{\gamma-1} \right)}&&\\
		&& \leq \E_{\Q} \left( \1_{\{ \underline{V}(w_n)\geq -a \}} \sum_{k=1}^{n}\ee^{-(\gamma-1)V(w_{k-1})} ( \sum_{u\in \Omega(w_k)} \ee^{-[V(u)-V(w_{k-1})]}   W_n^{(a)}(u)   )^{\gamma-1} \right)
		\\
		&&\leq  \E_{\Q} \left( \1_{\{ \underline{V}(w_n)\geq -a \}} \sum_{k=1}^{n}\ee^{-(\gamma-1)V(w_{k-1})} ( \sum_{u\in \Omega(w_k)} \ee^{-[V(u)-V(w_{k-1})]} \E\left( W_n^{(a)}(u) \big| \mathcal{G}_\infty \right) )^{\gamma-1} \right)
		\\
		&&=   \E_{\Q} \left( \1_{\{ \underline{V}(w_n)\geq -a \}}  \sum_{k=1}^{n}\ee^{-(\gamma-1)V(w_{k-1})} \sum_{u\in \Omega(w_k)} (\ee^{-[V(u)-V(w_{k-1})]} \P_{a+V(u)}\left( \underline{S}_{n-k}\geq 0  \right))^{\gamma-1} \right).
	\end{eqnarray*}
	By (\ref{chap3lldda2}) there exists $c>0$ such that
	$ 		\P_{a+V(u)}\left( \underline{S}_{n-k}\geq 0  \right)\leq c\frac{(1+a+V(u))}{\sqrt{n-k}}
	$, it stems that
	\begin{eqnarray*}
		\E_\Q((\mathtt{A})^{\gamma-1}) &\leq &   \E_{\Q} \left( \1_{\{ \underline{V}(w_n)\geq -a \}} \sum_{k=1}^{n}\frac{\ee^{-(\gamma-1)V(w_{k-1})}}{(n-k)^{\frac{\gamma-1}{2}}}    \Big(\sum_{u\in \Omega(w_k)} (1+a+V(u)) \ee^{-[V(u)-V(w_{k-1})]}  \Big)^{\gamma-1} \right)
	\end{eqnarray*}
	For any $u\in \omega(w_k)$, $(1+a+V(u))\leq (a+V(w_{k-1})_+) (1+ [V(u)-V(w_{k-1})]_+)$, so we can write
	\begin{eqnarray*}
		\E_\Q((\mathtt{A})^{\gamma-1}) &\leq&     c \sum_{k=1}^{n} \E_{\Q} \Big( \1_{\{ \underline{V}(w_n)\geq -a \}}  \frac{(a+V(w_{k-1}))^{\gamma-1}}{(n-k)^{\frac{\gamma-1}{2}}}   \ee^{-(\gamma-1)V(w_{k-1})} \qquad 
		\\
	&&\qquad	\times   \Big(  \sum_{u\in \Omega(w_k)} (1+[V(u)-V(w_{k-1})]_+)\ee^{-[V(u)-V(w_{k-1})]} \Big)^{\gamma-1}  \Big).
	\end{eqnarray*}
	Now, by using the branching property at time $k-1$ then (\ref{chap3lldda2}), we have (\ref{Pmoment}), we have
	\begin{eqnarray*}
		\E_\Q((\mathtt{A})^{\gamma-1})&\leq&     c \sum_{k=1}^{n} \E_{\Q} \Big( \1_{\{ \underline{V}(w_{k-1})\geq -a \}}  \frac{(a+V(w_{k-1}))^{\gamma-1}}{(n-k)^{\frac{\gamma-1}{2}}}   \ee^{-(\gamma-1)V(w_{k-1})} \qquad 
		\\
		&&\times  \E_\Q\Big( \frac{([a+V(w_{k-1})]\times (1+[V(w_1)]_+)}{(n-k)^{\frac{1}{2}}} \sum_{u\in \Omega(w_1)} (1+[V(u)]_+)\ee^{-V(u)]} \Big)^{\gamma-1}  \Big).
	\end{eqnarray*}
	By (\ref{Spiine}) and  (\ref{Pmoment}) we have 
	\begin{equation*}
		\E_\Q\Big((1+[V(w_1)]_+)\sum_{u\in \Omega(w_1)} (1+[V(u)]_+)\ee^{-V(u)]} \Big)^{\gamma-1}  \Big)\leq \E_\Q\Big( \sum_{|u|=1} (1+[V(u)]_+)\ee^{-V(u)} \Big)^{\gamma}  \Big)<+\infty
	\end{equation*}
	So by using again $(1+a+V(u))\leq (a+V(w_{k-1})_+) (1+ [V(u)-V(w_{k-1})]_+)$, we get
	\begin{eqnarray*}
		\E_\Q((\mathtt{A})^{\gamma-1}) &\leq &    c \sum_{k=1}^{n} \E_{\Q} \Big( \1_{\{ \underline{V}(w_{k-1})\geq -a \}}  \frac{(a+V(w_{k-1}))^{\gamma}}{(n-k)^{\frac{\gamma-1}{2}}}   \ee^{-(\gamma-1)V(w_{k-1})} \Big)
		\\
		&\leq & c \ee^{(\gamma-1) a} \E_{a} \left( \1_{\{ \underline{S}_{k-1}\geq 0\}}  \sum_{k=1}^{n}\frac{S_{k-1}^{\gamma}}{(n-k)^{\frac{\gamma}{2}}}   \ee^{-(\gamma-1)S_{k-1}}  \right).
	\end{eqnarray*}
	Finally by distinguish $k\leq \frac{n}{2}$ and  using Lemma B.2 in \cite{Aid11} or $k\geq \frac{n}{2}$ and using (\ref{chap3Aine}), it is not difficult to show that
	\begin{equation}
		n^{\gamma\over 2}\E_{\Q}((\mathtt{A})^{\gamma -1}) \leq c \ee^{(\gamma-1) a },
	\end{equation} 
	which concludes the proof of Lemma \ref{BRWPmoment} \hfill$\Box$
	\\

	\noindent{\it Proof of Proposition \ref{ConcentrationBRW}.}
	According to the Markov inequality then the branching property we have
	\begin{eqnarray*}
		\lefteqn{\P\left(  \left| \E\left( \sqrt{n+l}\, W_{n+l}^{(a)} \big|\mathcal{F}_l\right)-   (t+l)^{1/2}\, W_{n+l}^{(a)}\right|\geq \epsilon \right)}\\
		&&\epsilon^{-\gamma} \E\left(  \left| \sum_{|z|=l} \ee^{-V(z)}\1_{\{  \underline{V}(z)\geq -a \}}   \sqrt{n+l}\left[  \E\left( W_{n}^{(y+a)} \right) - W_n^{(y+a)}\right]_{\big| y=V(z)} \right|^\gamma \right).
	\end{eqnarray*}
	Applying inequality (\ref{russkov}), we deduce that there exists $c>0$ such that
	\begin{eqnarray*}
		\lefteqn{ \P\left(  \left| \E\left( \sqrt{n+l}\, W_{n+l}^{(a)} \big|\mathcal{F}_l\right)-   (t+l)^{1/2}\, W_{n+l}^{(a)}\right|\geq \epsilon \right)}\\		
		&& \leq c\epsilon^{-\gamma} \E\left( \sum_{|z|=l}\ee^{-\gamma V(z)} \1_{\{ \underline{V}(z)\geq -a  \}} \E\left(\left| \sqrt{n+l}\,  W^{(y+a)}_n\right|^\gamma \right)_{\big| y=V(z)}\right) 
		\\
		&&\leq c  \frac{ (n+l)^\frac{\gamma}{2} \ee^{(\gamma-1) a}}{n^{\frac{\gamma}{2}}\epsilon^\gamma  } \E\left( \sum_{|z|=l} \ee^{-V(z)} \1_{\{ \underline{V}(z)\geq -a\}} \right),
	\end{eqnarray*}
	where in the last inequality we have used (\ref{eqMomentdomination}). To conclude, we  apply the  many-to-one identity, (\ref{chap32.1}), and then (\ref{chap3lldda2}), to conclude that
	\begin{eqnarray*}
		\P\left(  \left| \E\left( \sqrt{n+l}\, W_{n+l}^{(a)} \big|\mathcal{F}_l\right)-   (t+l)^{1/2}\, W_{n+l}^{(a)}\right|\geq \epsilon \right) &\leq &  c  \frac{ (n+l)^\frac{\gamma}{2}}{n^{\frac{\gamma}{2}}\epsilon^\gamma \ee^{(\gamma-1) a} } \Q\left(   \underline{S}_n \geq -a  \right)
		\\
		& \leq &  \frac{c(n+l)^\frac{\gamma}{2}}{n^{\frac{\gamma}{2}}\epsilon^\gamma  } \frac{ (1+a) \ee^{(\gamma-1) a}}{\sqrt{l}},
	\end{eqnarray*}
	which ends the proof of  (\ref{SeneCondBRW}). \hfill$\Box$
	\\
	
	\noindent{\it Proof of Theorem \ref{BRWtheorem}.} By combining Proposition \ref{BRWlargenumber} and Proposition \ref{ConcentrationBRW}, for any $a>0$ we have
	\begin{equation}
		\lim_{ n\to\infty} \sqrt{n} W_n^{(a)}= \frac{1}{c_0} \sqrt{\frac{2}{\pi\sigma^2}} M_\infty'^{(a)},\qquad \text{in } \P \text{ probability}.
	\end{equation}
	The desired result now holds because  $\lim_{a\to \infty} M_\infty'^{(a)}=M_\infty'$ and \[
	 \lim_{a\to \infty} \lim_{ n\to\infty} \sqrt{n} W_n^{(a)}=\lim_{ n\to\infty} \sqrt{n} W_n^{(a)} ,\] 
	 almost surely in both cases.   \hfill$\Box$

	\subsection{Branching Brownian motion.}
	
	The case of branching Brownian motion is extremely similar to the case of branching random walks. For example, sampling a Branching Brownian motion on a lattice of times reveals an embedded branching random walk. This is not the approach one appeal to in order to  produce the desired results, however. One can mimic the proofs of the branching random walk in the previous section from first principles; the estimates and computations going through almost verbatim. We leave the details to the reader as an exercise. 
	
	\subsection{Gaussian multiplicative chaos}

	 The lognormal star-scale random measures were introduced in \cite{ARV11}. They are a very important class of random measures satisfying a continuous version of the celebrated Mandelbrot star equation. The authors of \cite{ARV11} characterize such measures by using the chaos theory of Kahane \cite{Kah85}.  Our aim is to reproduce a Seneta-Heyde result for these measures, originally due to \cite{DRSV12b}, using the generic approach highlighted in earlier sections, thereby offering a slimmer proof from the original. 
	
	Consider a family of centered stationary  Gaussian processes $(X_s(x))_{s\geq 0,\, x\in \r^{d}}$ $d\geq 1$, with covariances \nomenclature[e1]{$(X_s(x))_{s\geq 0,\, x\in \r^{d}}$}{: the log-correlated Gaussian field} \nomenclature[e2]{$\mathtt{k}$}{: the kernel function of $X_\cdot(\cdot)$}\nomenclature[e3]{$\mathtt{g}(\cdot)$}{$:=1-\mathtt{k}(\cdot)$}
	\begin{equation}
		\label{defX}
		\E[X_t(0)X_t(x)]=\int_0^t\mathtt{k}(\ee^u x)du,\quad \forall t>0,\, x\in \r^d.
	\end{equation}
	The kernel function $\mathtt{k}:\r^d\to \r$ is $\mathcal{C}^1$, satisfying $\mathtt{k}(0)=1$ and $\mathtt{k}(x)=0$ if {\bf $x\notin B(0,1):=\{x: |x|\leq 1\}$} ($|x|:= \underset{i\in [1,d]}{\max}|x_i|$). We also denote $\mathtt{g}(\cdot):=1- \mathtt{k}(\cdot)$ and introduce for any $t>0$,\nomenclature[e4]{$Y_t(x)$}{$:=X_t(x) -\sqrt{2d} t $}\nomenclature[c1]{$B(0,1)$}{$ :=\{x:\, \mid x \mid \leq 1\} $}
	\begin{equation}
		\label{defY}
		Y_t(x):=X_t(x)-\sqrt{2d}t.
	\end{equation}
	Let $\mathcal{B}(\r^d)$ the Borel on $\r^d$, and $\mathcal{B}_b(\r^d) $ its restriction to the bounded sets.  We introduce for $t>0$ and $\gamma>0$, the random measures $M_t'(dx)$ and $M_t^\gamma(dx)$ defined by:
	\begin{equation}
		\label{1.4} M_t'(A):=\int_{A}(-Y_t(x))\ee^{\sqrt{2d}Y_t(x)+d t}dx, \, \, \, M_t^{\sqrt{2d}}(A):=\int_{A} \ee^{\sqrt{2d}  Y_t(x)+ 2d t-dt}dx, \,\, \forall A\in \mathcal{B}_b(\r^d).
	\end{equation}
	are respectively the {\it critical} derivative and additive measure (when $A=[0,1]^d$ we will use the notations $M'_t$ and $M^{\sqrt{2d}}_t$). Kahane in \cite{Kah85} proved that 
	$$ \lim_{t\to \infty} M_t^{\sqrt{2d}}(A)=0,\qquad \text{a.s}\quad \forall A\in \mathcal{B}_b(\r), $$
	and recently Duplantier, Rhodes, Sheffield and Vargas \cite{DRSV12a} prove that there exists a random measure $M'$ such that 
	$$ \lim_{t\to \infty} M_t'(A)= M'(A),\qquad \text{a.s}\quad \forall A\in\mathcal{B}_b(\r^d).$$
	Our aim is to recover  the Seneta-Heyde norming first obtained in \cite{DRSV12b}, by using the same approach as for the fragmentation and branching random walk process.
	\begin{theorem}
		\label{SenetaGMC}
		The family $(\sqrt{t} M_t^{\sqrt{2d}})_t$ converges in probability as $t\to \infty$ towards a non trivial limit. More precisely,
		\begin{equation}
			\sqrt{t} M_t^{\sqrt{2d}}(A) \to \sqrt{\frac{2}{\pi}}M'(A),\qquad \text{in }  \P\text{ probability} \text{ as } t\to \infty.
		\end{equation}
	\end{theorem}

	\subsection{Proof of Theorem \ref{SenetaGMC}.}
	
	Let us introduce the truncated versions of the critical measures, $\forall A\in \mathcal{B}_b(\r^d)$, 
	\begin{equation}
		M_t^{\sqrt{2d}, (a)}(A):= \int_{A}  \ee^{\sqrt{2d} Y_t(x)+dt}\1_{\{ \sup_{s\leq t}Y_s(x) \leq a \}} dx,\quad M_t'^{(a)}(A):= \int_{A} (a-Y_t(x))\ee^{\sqrt{2d}Y_t(x)+dt} \1_{\{ \sup_{s\leq t}Y_s(x) \leq a \}} dx.
	\end{equation}
When $A=[0,1]^d$ we will use the notations $	M_t^{\sqrt{2d}, (a)} $ and $  M_t'^{(a)}$.

	\begin{Proposition}
		\label{GMClargenumber}
		For any $a,\, t,\, l>0,\, A\in \mathcal{B}_b(\r^d)$, we have
		\begin{equation}
			\lim_{t\to \infty} \sqrt{t} \E( M_{t+l}^{\sqrt{2d}, (a)}(A)  \big| \mathcal{F}_l)\overset{\text{a.s}}{=} {\sqrt{\frac{2}{\pi }}M_l'^{(a)}}(A).
		\end{equation}
	\end{Proposition}
	\noindent{\it Proof of Proposition \ref{GMClargenumber}.}
	By the Girsanov transformation we know that
	\begin{equation}
		\label{Sameway} 
		\E\left(  \ee^{\sqrt{2d} Y_t(x)+dt}\1_{\{ \sup_{s\leq t}Y_s(x) \leq a-y \}}   \right)= \P\left(  \sup_{s\leq t}B_s \leq a-y \right),\qquad \forall y\in \r. 
	\end{equation}
	Moreover by using the Fubini'Theorem, we get
	\begin{eqnarray*}
		\E( M_{t+l}^{\sqrt{2d}, (a)}(A)  \big| \mathcal{F}_l)&=& \sqrt{t} \int_{A}  \ee^{\sqrt{2d} Y_l(x)+dl}\1_{\{ \sup_{s\leq l}Y_s(x) \leq a \}} \E\left(  \ee^{\sqrt{2d} Y_t(x)+dt}\1_{\{ \sup_{s\leq t}Y_s(x) \leq a -y \}}   \right)_{\big|y=Y_l(x)} dx
		\\
		&= & \int_{A}  \ee^{\sqrt{2d} Y_l(x)+dl}\1_{\{ \sup_{s\leq l}Y_s(x) \leq a \}}\sqrt{t} \P\left(  \sup_{s\leq t}B_s \leq a -y  \right)_{\big|y=Y_l(x)} dx
		\\
		&\underset{t\to \infty}{\longrightarrow}&  \sqrt{\frac{2}{\pi}} \int_{A} (a-Y_l(x))  \ee^{\sqrt{2d} Y_l(x)+dl}\1_{\{ \sup_{s\leq l}Y_s(x) \leq a \}}   dx =    \sqrt{\frac{2}{\pi}}  M_l'^{(a)}(A),
	\end{eqnarray*}
	which proves the Proposition \ref{GMClargenumber}.\hfill$\Box$
	
	\begin{Proposition}
		\label{ConcentrationGMC} Let $A$ a Borel set of $[0,1]^d$. There exists $c>0$, $\gamma\in (1,2)$ such that:  For any $\epsilon >0,\, l,\, t\geq 0$ we have
		\begin{equation}
			\label{SeneCondGMC}
			\P\left(  \Big| \E\left( (t+l)^{1/2}\, M_{t+l}^{\sqrt{2d}, (a)}(A) \big|\mathcal{F}_l\right)-   (t+l)^{1/2}\, M_{t+l}^{\sqrt{2d}, (a)}(A)\Big|\geq \epsilon \right)\leq \frac{c(t+l)^\frac{\gamma}{2}}{t^{\frac{\gamma}{2}}\epsilon^\gamma \sqrt{l}}(1+a)\ee^{\sqrt{2d}  a}.
		\end{equation}
	\end{Proposition}

	\begin{Lemma}
		\label{GMCPmoment}
		Let $O$ a open set of $[0,1]^d$, and $\rho,\, l$ two continuous function. There exists $c>0$ such that for any $t\geq 1$ large enough, we have
		\begin{equation}
			\E\left(  \left( \sqrt{t} \int_{O}\ee^{\sqrt{2d}[ Y_t(x)+l(x)] +d t}\1_{\{ \sup_{s\leq t} Y_s(x)\leq \rho(x) \}} dx \right)^2 \right)\leq c \int_{x\in O}\ee^{\sqrt{2d} [\rho(x)+2l(x)]}dx.
		\end{equation}
	\end{Lemma}
	\noindent{\it Proof of Lemma \ref{GMCPmoment}.}  Let $K$ large enough such that for any $x,y\in [0,1]^d$, $|x-z|\leq K^{-1}$ implies $|\rho(x)-\rho(z)|+|l(x)-l(z)|\leq 1$.
	\begin{eqnarray*}
		&& t \int_{O}\int_{O} \E\left(  \ee^{\sqrt{2d}[Y_t(z)+l(z)] +dt}\1_{\{  \sup_{s\leq t} Y_s(z)\leq \rho(z) \}}        \ee^{\sqrt{2d}[Y_t(x)+l(x)] +dt}\1_{\{  \sup_{s\leq t} Y_s(x)\leq \rho(x) \}}  \right)  dx   dz
		\\
		&&\leq t[(\mathtt{S}) +(\mathtt{A})],
	\end{eqnarray*}
	with 
	\begin{eqnarray}
		\label{spineplusA}	&&(\mathtt{S}):=  \int_{O}\int_{x\in [0,1]^d,\, |x-z|\leq \ee^{-t} }   \E^{*}(x,z)  dx   dz,\qquad (\mathtt{A}):=  \int_{O}\int_{ x\in [0,1]^d,\, |x-z|\geq \ee^{-t} } \E^{*}(x,z) dx   dz
	\end{eqnarray}
	and
	\begin{equation}
		\label{ardu} \E^{*}(x,z):= \E\left(  \ee^{\sqrt{2d}[Y_t(z)+l(z)+\sqrt{2d}Y_t(x)+l(x)] +2dt}\1_{\{  \sup_{s\leq t} Y_s(z)\leq \rho(z) ,\,   \sup_{s\leq t} Y_s(x)\leq \rho(x) \}}  \right).
	\end{equation}
	We first bound $(\mathtt{A})$, the ``additive" term, it requires to estimates $ \label{ardu} \E^{*}(x,y)$ for $|x-y|\geq \ee^{-t}$. Let introduce some notations. For any $x,z\in [0,1]^d$ let denote $\lambda_{x,y}:= -\ln {|x-y|}$. For $t\geq 0$,  let $\mathcal{F}_t$ be the sigma-field generated by the Gaussian field $(Y_u(x))_{u\leq t,\, x\in \r^d}$.   According to (\ref{defX}), it is plain to check that for any  $|x-z|\geq \ee^{-t}$, the two processes 
	$$(Y_s(z)-Y_{\lambda_{x,y}}(z))_{s\in [\lambda_{x,z} ,t]},\quad \text{and}\quad  (Y_s(x)-Y_{\lambda_{x,z}}(x))_{s\in [\lambda_{x,z},t]},$$ 
	are independent. We take advantage of this independence by first conditioning the expectation that defines $\E^{*}(x,z)$  on   the sigma field $\mathcal{F}_{\lambda_{x,z}}$, by using in addition the Girsanov' identity in the same way as (\ref{Sameway}), it yields
	\begin{eqnarray*}
		\E^{*}(x,y)\leq  c\ee^{\sqrt{2d}[l(x)+l(z)]} \E\Big(  \ee^{\sqrt{2d}[Y_{\lambda_{x,z}}(z) +Y_{ \lambda_{x,z} }(x) +2d \lambda_{x,z} ]}\1_{\{  \sup_{s\leq   \lambda_{x,z}} Y_s(z)\leq \rho(z) ,\,  \sup_{s\leq  \lambda_{x,z}} Y_(x)\leq \rho(x) \}}\times 
		\\
		\P_{ Y_{ \lambda_{x,z}}(z)-\rho(z) }\left( \sup_{s\leq t-  \lambda_{x,z}  }B_s\leq 0 \right)   \P_{Y_{ \lambda_{x,z} }(x) -\rho(x)}\left( \sup_{s\leq t- \lambda_{x,z}}B_s\leq 0\right) \Big). 
	\end{eqnarray*}
	Using now (\ref{chap3lldda2}), we get
	
	\begin{eqnarray*}
		\E^{*}(x,z)\leq 	c\frac{\ee^{\sqrt{2d}[l(x)+l(z)]}}{t-\lambda_{x,z}  }   \E\left(  (\rho(z) -Y_{  \lambda_{x,z} }(z) +1)(\rho(x) -Y_{ \lambda_{x,z} }(x) +1)\ee^{\sqrt{2d}[Y_{ \lambda_{x,z} }(z)+ Y_{ \lambda_{x,z} }(x)] +2d  \lambda_{x,z}  } \right.
		\\
		\left. \times\1_{\{  \sup_{s\leq   \lambda_{x,z}  } Y_s(z)\leq \rho(z) ,\,  \sup_{s\leq   \lambda_{x,z} } Y_(x)\leq \rho(x) \}}   \right).
	\end{eqnarray*}
	Moreover by using  the Cauchy-Schwartz inequality, it appears two terms of this type
	$$ \ee^{\sqrt{2d}l(x)}  \E\left( (\rho(x) -Y_{ \lambda_{x,z} }(x) +1)^2 \ee^{2\sqrt{2d} Y_{ \lambda_{x,z} }(x) + 2d  \lambda_{x,z}  } \1_{\{  \sup_{s\leq   \lambda_{x,z} } Y_(x)\leq \rho(x) \}}    \right)^{\frac{1}{2}}. $$
	Hence, using the Girsanov' transformation we can affirm that 
	$$\E^{*}(x,y)\leq c \frac{F(x,\lambda_{x,z}) F(z,\lambda_{x,z})}{t-\lambda_{x,y}} . $$ 
	with
	\begin{eqnarray}
		F(z,r)&:=& \ee^{\sqrt{2d}  l(z)+ \frac{\sqrt{2d}}{2}\rho(z)} \ee^{\frac{d}{2}r} \E_{-\rho(z)-1}\left(  (B_{r })^2   \ee^{\sqrt{2d}B_{r}(z)  } \1_{\{  \sup_{s\leq  r} B_s\leq  0 \}}    \right)^\frac{1}{2}.
	\end{eqnarray}
	Going back to the definition of $(\mathtt{A})$ in (\ref{spineplusA}), we get that
	\begin{eqnarray*}
		(\mathtt{A}) &\leq&  \left(\int_{O} \int_{O}\1_{\{ |x-z|\geq \ee^{-t} \}}   \frac{F(x,\lambda_{x,z}) F(z,\lambda_{x,z})}{t-\lambda_{x,y}}   dxdz\right).
	\end{eqnarray*}
	Noticing the symmetry between $x$ and $z$ and applying the Cauchy-Schwartz inequality we obtain
	\begin{eqnarray*}
		(\mathtt{A})	&\leq & c\left(\int_{O} \int_{O}\1_{\{ |x-z|\geq \ee^{-t} \}}   \frac{F(x,\lambda_{x,z})^2}{t-\lambda_{x,y}}   dxdz\right)
		\\
		&\leq & c'  \int_{O} \ee^{2\sqrt{2d} l(z)+\sqrt{2d} \rho(z) } \int_{O}\1_{\{ |x-z|\geq \ee^{-t} \}}  \frac{|x-y|^{-d}}{t-\lambda_{x,z}} \E_{-\rho(z)-1}\left(  (B_{r })^2   \ee^{\sqrt{2d}B_{r}(z)  } \1_{\{  \sup_{s\leq  r} B_s\leq  0 \}}    \right)  dx dz.
	\end{eqnarray*}
	Finally by decomposing $O$, as the union of the annulus  $(\mathtt{C}_k(x))_{ 0\leq k \leq  \lfloor t\rfloor +1}$, which are centered in $x$ and with radius $\ee^{-k}$ and $\ee^{-(k-1)}$, then using (\ref{eqPourRW2}), it yields
	\begin{eqnarray}
		(\mathtt{A})	&\leq & \frac{c}{t} \int_{O}  \ee^{2\sqrt{2d} l(z)+\sqrt{2d} \rho(z) } dz.
	\end{eqnarray}
	Concerning $(\mathtt{S})$, the ``spine term", when $t$ is large enough,  for any $x,y\in [0,1]^d$ such that $|x-z|\leq \ee^{-t}\leq K^{-1}$ we have $|\rho(z)-\rho(x)|+|l(x)-l(z)|\leq 1$ and get
	\begin{eqnarray*}
		\E^{*}(x,y) &\leq& \E\left(  \ee^{2\sqrt{2d}[Y_t(z) +l(z)] +2dt}\1_{\{  \sup_{s\leq t} Y_s(z)\leq \rho(z)  \}}  \right)^\frac{1}{2} \E\left(  \ee^{2\sqrt{2d}[ Y_t(x)+l(x)]  +2dt}\1_{\{    \sup_{s\leq t} Y_s(x)\leq \rho(x) \}}  \right)^\frac{1}{2}
		\\
		&\leq& c'\ee^{dt} \ee^{2\sqrt{2d} l(z)}\E\left( \ee^{\sqrt{2d}B_t}\1_{\{  \sup_{s\leq t} B_s\leq \rho(z)+1 \}}   \right).
	\end{eqnarray*}
	And finally we have
	\begin{eqnarray*}
		(\mathtt{S})&\leq &\int_{O}\ee^{2\sqrt{2d} l(z)+ \sqrt{2d}\rho(z)} \int_{O} \ee^{dt} \1_{\{|x-z|\leq \ee^{-t}\}}   \E_{-\rho(z)-1}\left( \ee^{\sqrt{2d}B_t}\1_{\{  \sup_{s\leq t} B_s\leq 0 \}}   \right)  dx  dz
		\\
		&\leq & \frac{c}{t^\frac{3}{2}} \int_{O}  \ee^{2\sqrt{2d} l(z)+\sqrt{2d} \rho(z) } dz,
	\end{eqnarray*}
	which finishes the proof of Lemma \ref{GMCPmoment} \hfill$\Box$
	\\

	\noindent{\it Proof of Proposition \ref{ConcentrationGMC}.} By (\ref{defX}) we can easily check that for any $l>0$, the process
	\begin{equation}
		\label{increment}( Y_s^{(l)}(x))_{s\geq 0,\, x\in \r^d}:= ( Y_{s+l}(x)-Y_l(x))_{s\geq 0,\, x\in \r^d},
	\end{equation}
is independent of $( Y_s(x))_{s\leq l,\, x\in \r^d}  $ and has the same law as $( Y_s(x\ee^l))_{s\geq 0,\, x\in \r^d}$. Moreover for any  $|x-z|\geq \ee^{-l}$, the two processes $  ( Y_s(x))_{s\leq l,\, x\in \r^d},\quad \text{and}\quad   ( Y_s(z))_{s\leq l,\, x\in \r^d},$ are independent. To take profit of these two properties of the field $Y$, following \cite{Mad12} we shall introduce a buffer zone. 

Let $l>0$, we will divide the cube $[0,\ee^{l}]^d$ into several smaller cube of size $R>0$, all of these smaller cube being at distance greater than $1$ from each other. We assume that $R,l$ are such that
$$
m:= \frac{ \ee^{l} +1  }{R+1}\in \N^*.$$
The integer $m$ stands for the number of small squares of size $R$ that one meets along an edge of the  cube $[0,\ee^l]^d$. The basis of each small square will be indexed with a $d$-uplet
	$${\bf i}= (i_1,...,i_d)  \in\{1,...,m\}^d    .$$ The basis of the square $D_{{\bf i}}$ is then located at 
	$$a_{{\bf i}}:= (R+1)\left((i_1-1),...,(i_d-1)\right)\in [0,e^{l}]^d  .$$
	in such a way that  
	$$D_{{\bf i}}:= a_{\bf i}+ [0,R]^d.$$
	One may observe that all the squares $D_{{\bf i}}$ are separated from each other by a fishnet shaped buffer zone, which is precisely
	$${\rm BZ}_{R,l}:=[0,e^{l}]^d\setminus \bigcup_{{\bf i}\in \{1,...,m\}^d} D_{{\bf i}}.$$

	The terminology "buffer zone" is used because this is the minimal area needed to make sure that the values taken by the process $Y_.^{(l)}$ inside each $D_{\bf i}$ are independent of its values on all other $D_{\bf j}$ for ${\bf j}\neq {\bf i}$. Now let us define
	\begin{equation}
		\tilde{M}_{t+l}^{\sqrt{2d},(a)}(A):= \sum_{{\bf i}}  M_{t+l}^{\sqrt{2d},(a)}(A\cap (D_{\bf i}. \ee^{-l}))=  \sum_{{\bf i}} \int_{A\cap (D_i .\ee^{-l})} \ee^{\sqrt{2d} Y_{t+l}(x)+d(t+l)}\1_{\{  \sup_{s\leq t+l} Y_s(x)\leq a    \}},
	\end{equation}
	where $D_{\bf i}:= \{ (x_1,...,x_d)\in \r^d,\, (x_1\ee^{l},...,x_d\ee^{l})\in D_{\bf i} \}$. We are now ready to tackle the proof of the Proposition \ref{ConcentrationGMC}. By the Markov inequality,
\begin{equation}
	\label{hoss}
		\P\left(  \left| \E\left( (t+l)^{1/2}\, M_{t+l}^{\sqrt{2d},(a)} \big|\mathcal{F}_l\right)-   (t+l)^{1/2}\, M_{t+l}^{\sqrt{2d}, (a)}\right|\geq \epsilon \right)\leq 4\epsilon^{-2}\E\left(  M_{t+l}^{\sqrt{2d},(a)}(A)-  \tilde{M}_{t+l}^{\sqrt{2d},(a)}(A) \right) +\P^*
	\end{equation}
	with 
	$$\P^*: = \P\left(  \left| \E\left( (t+l)^{1/2}\, \tilde{M}_{t+l}^{\sqrt{2d},(a)}(A) \big|\mathcal{F}_l\right)-   (t+l)^{1/2}\, \tilde{M}_{t+l}^{\sqrt{2d}, (a)}(A)\right|\geq \frac{\epsilon}{2} \right).$$
By taking the conditional expectation with respect to $\mathcal{F}_l$ then by using the inequality (\ref{russkov}), we have
	\begin{eqnarray*}
	\P^*	&\leq &\P\left(  \left| \E\left( (t+l)^{1/2}\, \tilde{M}_{t+l}^{\sqrt{2d},(a)}(A) \big|\mathcal{F}_l\right)-   (t+l)^{1/2}\, \tilde{M}_{t+l}^{\sqrt{2d}, (a)}(A)\right|\geq \frac{\epsilon}{2} \right) 
		\\
		&\leq&  \epsilon^{-2} \sum_{{\bf i} } \E\left( \E\left( \left(   M_{t+l}^{\sqrt{2d},(a)}(A \cap (D_{\bf i}. \ee^{-l}))  \right)^2 \big|\mathcal{F}_l\right)\right) 
		\\
	&	\leq& \frac{(t+l)^{1/2}}{\sqrt{t}} \epsilon^{-2} \sum_{{\bf i}} \E\left( \E\left( \left( \int_{D_{\bf i}\cdot  \ee^{-l} }  \sqrt{t} \ee^{\sqrt{2d} [Y_{t}(x)+\lambda(x) ]+dt}\1_{\{  \sup_{s\leq t} Y_s(x)\leq a  +\rho(x)   \}}    \right)^2 \right)_{ \Big|  \rho(\cdot)=-\lambda(\cdot)= -Y_l(\cdot\ee^{-l})  }\right).
	\end{eqnarray*}
	Then by Lemma \ref{GMCPmoment} (applied with $O =A\cap D_i\cap\{x\in [0,1]^d,\, \sup_{s\leq l} Y_s(x)\leq a \}$), we obtain
	\begin{eqnarray}
		\nonumber \P^* &\leq & \sum_{{\bf i} } \E\left(\int_{D_i}  \ee^{\sqrt{2d}Y_l(x \ee^{-l})} \1_{\{  \sup_{\leq l} Y_s(x)\leq a\}}dx \right) \ee^{\sqrt{2d} a}
		\\
		\nonumber  &\leq &  \ee^{\sqrt{2d} a} \frac{(t+l)^{1/2}}{\sqrt{t}}\sum_{{\bf i}}\lambda(D_{\bf i}\ee^{-l}) \P\left( \sup_{s\leq l} B_s\leq a  \right)
		 \\
	\label{corepf}	 &\leq&  \ee^{\sqrt{2d} a}(1+a) \frac{(t+l)^{1/2}}{\sqrt{t}} \frac{c\lambda([0,1]^d)}{\sqrt{l}}.
	\end{eqnarray}
where, in the second inequality,  we used the change of variables $x=\ee^{-l}y$ and in the last  inequality we used (\ref{chap3lldda2}). Otherwise we have
	\begin{eqnarray}
		\nonumber 4\epsilon^{-2}\E\left(  M_{t+l}^{\sqrt{2d},(a)}(A)-  \tilde{M}_{t+l}^{\sqrt{2d},(a)}(A) \right)& = & \epsilon^{-2}\int_{BZ_{R,l}}\E\left(  \ee^{\sqrt{2d} Y_{t+l}(x)+d(t+l)}\1_{\{\sup_{s\leq t} Y_s(t)\leq a \}} \right) 
		\\
		\label{Buffer}&\leq& \epsilon^{-2} (1+a) \lambda(BZ_{R,l}) .
	\end{eqnarray}
	which goes to $0$ when $R$ goes to infinity. By combining (\ref{corepf}) and (\ref{Buffer}) we obtain the Proposition \ref{ConcentrationGMC}. \hfill$\Box$
	\\

	\noindent{\it Proof of Theorem \ref{SenetaGMC}.}
	Let $A\in \mathcal{B}([0,1]^d)$, by combining Proposition \ref{GMClargenumber} and \ref{ConcentrationGMC}, we have obtained that for any fixed $a>0$,
	\begin{equation}
	\label{oulou}
		\lim_{t\to\infty} \sqrt{t} M_t^{\sqrt{2d},(a)}(A)= \sqrt{\frac{2}{\pi}}M'^{(a)}(A),\qquad \text{in } \P \text{ probability}.
	\end{equation}
	Moreover it is known (see \cite{DRSV12b}) that a.s $\sup_{s\geq 0,\, x\in [0,1]^d}Y_s(x) <\infty$. We deduce that uniformly in $t\in \r$, 
	\begin{equation}
	\label{boulou}
		\lim_{a\to \infty} M_t'^{(a)}(A)=M_t'(A),\qquad \P \text{ a.s.}
	\end{equation}
Both (\ref{oulou}) and (\ref{boulou}) conclude the proof of the Theorem \ref{SenetaGMC}.\hfill$\Box$

\appendix

\section*{Appendix: Technical box}
\addcontentsline{toc}{section}{Appendix: Technical box}
\renewcommand{\thesubsection}{\Alph{subsection}}

In this section, let $Z: = ({Z}_s)_{s\geq 0}$ under $\Pp$ be a spectrally positive  L\'evy process. As usual, we shall write 
 $\Pp_a$ for the law of $Z$ conditional ${Z}_0=a$. We make the additional assumptions that $Z$ has with mean $0$ and variance equal to $\mathbb E({Z}_1^2)=\sigma^{2}$. Note that  these two assumption means that $Z$ is also a square integrable martingale. 
Moreover, there exists an  $\epsilon_0>0$ such that $\mathbb{E}(\ee^{\epsilon_0 |\xi_t|})<\infty. $ Finally we shall write  $\psi$ for the Laplace exponent of ${Z}$. That is to say, 
$\mathbb{E}\left( \ee^{-q {Z}_t }\right) =\exp( \psi(q) t),\, \forall q\geq 0,\, t\geq 0 $. Furthermore for $p,q\in \N$, we will use the notation  $[|p ,q|]:= \{k\in \N,\, q\leq k\leq p  \}$.

 \renewcommand{\theequation}{A.\arabic{equation}}
  \renewcommand{\thetheorem}{A.\arabic{theorem}}

  \setcounter{equation}{0}  
    \setcounter{theorem}{0}  

 \subsection{ L\'evy processes with no negative jumps conditioned to stay positive}
%

%

 %
  
 We are interested in understanding the asymptotic behaviour of spectrally positive L\'evy processes which survive crossing into the negative half-line in the long term and functional limits thereof. Recall the notation 
 \[
 \zeta^a = \inf\{s>0: Z_s<-a\}, \qquad a\in\mathbb{R}.
 \]

 \begin{Proposition}
 \label{Equivinfpos}For any $a>0$, we have
 \begin{equation}
 \label{Tool1}
 \lim_{t\to\infty} \sqrt{t}\Pp\left(\zeta^a>t \right)=  \sqrt{\frac{2}{\sigma^2 \pi}}\times a.
 \end{equation}
 \end{Proposition}
 \noindent{\it Proof of Proposition \ref{Equivinfpos}.} As $({Z}_s)_{s\geq 0}$ is a spectrally positive L\'evy process, $(\zeta^a)_{a\geq 0}$ is a subordinator (see VII.1 Theorem 1 in \cite{Ber98}) whose  Laplace exponent, $\phi$,  satisfies
 \begin{equation}
 \label{inverse}
 \psi(\phi(q))=q,\qquad \forall q>0.
 \end{equation}
As $\psi(0)=\psi'(0)=0$ and $\sigma^2 =\psi''(0)$, it follows that $\phi(q) {=} \sigma^2\frac{q^{2}}{2}+ o(q^{2})$ as $q\to0$, and thus
 \begin{equation}
 \phi(q)\underset{q\to 0}{=} \sqrt{\frac{2}{\sigma^2}} q^{\frac{1}{2}}+ o(q^{\frac{1}{2}}).
 \end{equation}
 Otherwise 
 \begin{eqnarray}
 1-\mathbb{E}(\ee^{-q\tau_a})= 1-\exp(-\phi(q)a) \underset{q\to 0}{=} a \sqrt{\frac{2}{\sigma^2}} q^{\frac{1}{2}}+ o(q^{\frac{1}{2}})
 \end{eqnarray}
and
\begin{eqnarray}
1-\mathbb{E}(\ee^{-q\tau_a})=q \int_0^\infty \ee^{-q u} \Pp(u\leq \tau_a) du.
\end{eqnarray}
So according to the Tauberian Theorem on pp. 10  of \cite{Ber98}, we have
$$ \int_{0}^x \Pp(u\leq \tau_a) du \underset{x\to\infty}{\sim} a \sqrt{\frac{2}{\sigma^2}}  \frac{x^{\frac{1}{2}}}{\Gamma(1+\frac{1}{2})}. $$
Then by the monotone density (e.g. on pp. 10 of \cite{Ber98}), we  can  affirm that
\begin{equation}
\Pp\left(  u\leq \tau_a  \right) \underset{u\to\infty}{\sim} \frac{au^{-\frac{1}{2}  }2}{\sqrt{2\sigma^2} \sqrt{\pi}}=\sqrt{\frac{2}{\sigma^2 \pi}}a u^{-\frac{1}{2}} .
\end{equation}
This concludes the proof of the inequality (\ref{Tool1}).
 \hfill$\Box$
\\
\begin{Proposition}[Donsker Theorem for  L\'evy process]
\label{Donsker22}
For any $F\in \mathcal{C}_b(\mathcal{D},[0,1])$ we have
\begin{equation}
\lim_{t\to \infty}\mathbb{E}_a\left( F\left( {t^{-1/2}}({Z}_{st})_{s\in [0,1]}  \right)   \right)= \mathbb{E}\left( F(\sigma (B_s)_{s\in [0,1]} \right).
\end{equation}
\end{Proposition}
\noindent{\it Proof of Proposition \ref{Donsker22}.} According to Theorem 15.5 in \cite{Bil99} it suffices to prove the finite-dimensional convergence of $({Z}_{st})_{s\in [0,1]}$ as well as verifying  the following conditions, which are equivalent to the tightness of $({Z}_{st})_{s\in [0,1]}$:
\begin{description}
\item[(i)]{\it For each positive $\eta$, there exists an $A$ such that 
\begin{equation}
\Pp_a\left( \frac{|{Z}_{0}|}{\sqrt{t}} \geq A    \right)\leq \eta,\qquad t\geq 0.
\end{equation}
}
\item[(ii)]{\it For each positive $\epsilon$ and $\eta$, there exist $\delta\in (0,1)$, and $t_0>0$ such that 
\begin{equation}
\Pp_a\left(  \sup_{s_1,s_2\in [0,1],\, |s_1-s_2|\leq \delta}{t^{-1/2}}\left|{Z}_{s_1 t}-{Z}_{s_2 t} \right|\geq \epsilon \right)\leq \eta,\qquad \forall t\geq t_0.
\end{equation}}

\end{description}
As $({Z}_{st})_{s\in [0,1]}$ is a process with independent and identically distributed increments, the finite dimensional convergence follows from the Lindeberg--L\'evy central limit theorem
Condition (i) is trivially satisfied, thus is suffices to prove (ii). Observe that
\begin{eqnarray*}
\Pp_a\left(  \sup_{s_1,s_2\in [0,1],\, |s_1-s_2|\leq \delta}{t^{-1/2}}\left|{Z}_{s_1 t}-{Z}_{s_2 t} \right|\geq \epsilon \right) &\leq & \sum_{i=0}^{\delta^{-1}} \Pp_a\left(  \sup_{ s \leq \delta t} \left|{Z}_{i\delta t+s}-{Z}_{i\delta t} \right|\geq \frac{\epsilon}{2} \sqrt{t} \right)
\\
&\leq & \delta^{-1} \Pp\left(  \sup_{s\leq \delta t} |{Z}_s|\geq \frac{\epsilon}{2} \sqrt{t}     \right).
\end{eqnarray*}
Moreover, since our assumptions mean that ${Z}$ is a square integrable martingale,  by Doob's martingale inequality, we have
\begin{eqnarray*}
\Pp_a\left(  \sup_{s_1,s_2\in [0,1],\, |s_1-s_2|\leq \delta}{t^{-1/2}}\left|{Z}_{s_1 t}-{Z}_{s_2 t} \right|\geq \epsilon \right) &\leq & \delta^{-1} \frac{2^4}{t^2\epsilon^4}\mathbb E\left( {Z}_{\delta t}^4 \right)
\\
&\leq & c\delta^{-1} 2^4 \frac{\delta^2}{\epsilon^4}\to 0,
\end{eqnarray*}
as  $\delta$ goes to $0$.\hfill$\Box$


\begin{Proposition} 
\label{Lemapart} 
Fix $a>0$.
\begin{description}
\item[(i)] We have,
  \begin{equation}
  \label{bpart} \lim_{K\to \infty} \lim_{t\to \infty} \mathbb{E}_a\left(   \frac{| {Z}_t|}{\sqrt{t}}  \1_{\{   \frac{| {Z}_t|}{\sqrt{t}}\geq K  \}}  \big|  \inf_{s\leq t} {Z}_s \geq 0 \right)=0.
  \end{equation}
  
\item[(ii)] For any $z\geq 0$, $l>0$ and $g\in \mathcal{D}([0,l])$, $F\in \mathcal{C}_b(\mathcal{D},[0,1])$, we have
 \begin{equation}
\label{apart}  \lim_{t\to\infty} \mathbb{E}_{a}\left(    F\left( {t^{-1/2}}\Big(z+   g_{st}\1_{\{ s\leq \frac{l}{t} \}}   +{Z}_{st}    \1_{\{ s\geq \frac{l}{t} \}}\Big)_{s\in [0,1]} \right)\big| \inf_{s\leq t} {Z}_s\geq 0 \right)
 =   \mathbb{E}\big(F((\sigma m_s)_{s\in [0,1]})\big),
 \end{equation} 
 with $(m_s)_{s\in [0,1]}$ a Brownian meander.
 \end{description}
\end{Proposition}
Recall that the under the probability 
\begin{equation}
\label{BessMeander}
\sqrt{\frac{2}{\pi}} m_1\cdot \mathbb{P},
\end{equation}
the process $(m_s)_{s\in [0,1]}$ has the same law as the three dimensional Bessel process. 

\begin{Remark}\rm An equivalent result is true for the random walk, it was first proved by Iglehart \cite{Igl74} then by Bolthausen \cite{Bol78}. Although it is intuitively clear that the result should hold in our setting, the extension to our process $({Z}_s)_{s\geq 0}$ is not straightforward. 
\end{Remark}

\noindent{\it Proof of Proposition \ref{Lemapart}}. The convergence (\ref{apart}) is the main part of the Proposition \ref{Lemapart}. First we shall prove (\ref{bpart}). For any $t,k>0$ observe that
\begin{eqnarray*}
\Pp_a\left(    \frac{| {Z}_t|}{\sqrt{t}}\geq k  ,  \inf_{s\leq t} {Z}_s \geq 0 \right) &\leq& \Pp_a\left(   \frac{|{Z}_{\lfloor t\rfloor}| +| {Z}_t - {Z}_{\lfloor t\rfloor}|}{\sqrt{t}}\geq k  , \min_{i \in [|1,\lfloor t\rfloor|]} {Z}_i \geq 0 \right)
\\
&&\leq \Pp_a\left(   \frac{|{Z}_{\lfloor t\rfloor}|}{\sqrt{t}}\geq \frac{k}{2}  , \min_{i \in [|1,\lfloor t\rfloor|]} {Z}_i \geq 0 \right)\\
&&+ \Pp_a\left( \min_{i \in [|1,\lfloor t\rfloor|]} {Z}_i \geq 0 \right)\Pp\left( | {Z}_t - {Z}_{\lfloor t\rfloor}|\geq \frac{k}{2}{\sqrt{t}} \right).
\end{eqnarray*}
By (\ref{Tool1}) and the assumed exponential moments of ${Z}$, for $t$ large enough, the second term is smaller than $\frac{c}{\sqrt{t}}\times \ee^{-\epsilon \frac{k}{2}\sqrt{t}}$. The process $({Z}_{i})_{i\in \N}$ is a standard random walk with $0$ mean. Thanks to the assumptions on ${Z}$, we know that $\mathbb{E}({Z}_1^2)<\infty$ and $\mathbb{E}(\ee^{\epsilon {Z}_1})<\infty$, so by (5.14) in \cite{Mad12} we deduce that the second term is smaller than $c\ee^{-c'k }$. We obtain that 
\begin{equation}
\label{eqaudessus} \forall t,k>0,\qquad \Pp_a\left(    \frac{| {Z}_t|}{\sqrt{t}}\geq k  ,  \inf_{s\leq t} {Z}_s \geq 0 \right)\leq \frac{c}{\sqrt{t}} \ee^{-c'k}.
\end{equation}
Now to prove (\ref{bpart}) it suffices to observe that
\begin{eqnarray*}
\mathbb{E}_a\left(   \frac{| {Z}_t|}{\sqrt{t}}  \1_{\{   \frac{| {Z}_t|}{\sqrt{t}}\geq K  \}}  \big|  \inf_{s\leq t} {Z}_s \geq 0 \right) &\leq& \sqrt{t} \sum_{k\geq K} (k+1) \Pp_a( |{Z}_t|\geq k,\, \inf_{s\leq t} {Z}_s\geq 0)
\\
&\leq &\sum_{k\geq K} \ee^{-c'k} \underset{K\to \infty}{\to } 0.
\end{eqnarray*}

Now we will prove (\ref{apart}). For any $z\geq 0$, $g\in \mathcal{D}([0,l])$, ${t^{-1/2}}(z + g)\to 0$ when $t$ goes to infinity. Moreover, as $F$ is a continuous function, it is sufficient to prove that $\forall a\geq 0$, $F\in \mathcal{C}(D,\r_+)$, 
\begin{equation}
\label{apartsimpler} \lim_{t\to\infty} \mathbb{E}_a\left(   F\{  {t^{-1/2}} ({Z}_{st})_{s\in[0,1]} \}\Big| \inf_{s\leq t} {Z}_s\geq 0 \right)= \mathbb{E}( F(\sigma m_s)_{s\in [0,1]}).
\end{equation}
For any $t,\, a>0$ we define
\begin{equation*}
\mathbb{P}_a^{(t)}(A):= \frac{ \mathbb{P}_a( ({Z}_{st})_{s\in [0,1]}\in A, \inf_{s\leq t} {Z}_s\geq 0  )}{\mathbb{P}_a(\inf_{s\leq t} {Z}_s\geq 0)},\qquad A\in \mathcal{B}(D).
\end{equation*}
By using the formalism of Billingsley \cite{Bil99}, we see that (\ref{apartsimpler}) is equivalent to the weak convergence of the sequence of probability measure $(\Pp_a^{(t)}(\cdot) )_{t\geq 0}$ to $\mathbb{W}^+$ the law of the Brownian meander.\\

 For any $(t_1,...,t_k)\in [0,1]^k$ let
$$ \pi_{t_1,...,t_k}: (f_s)_{s\in [0,1]}\mapsto (f_{t_1},...,f_{t_k}).  $$
It is well known (see Theorem 13.1 of \cite{Bil99}) that it suffices to show the tightness of  the sequence $(\Pp_a^{(t)}(\cdot))$ and its finite dimensional convergence to $\mathbb{W}^+$, i.e $\forall t_1,...,t_d\in [0,1]$, $$\Pp_a^{(t)}(\pi^{-1}_{t_1,...,t_d})\Rightarrow \mathbb{W}^+(\pi^{-1}_{t_1,...,t_d}).$$
\begin{Lemma}
\label{Tightness}
The sequence $(\Pp_a^{(t)}(\cdot))$ is tight.
\end{Lemma}
\noindent{\it Proof of Lemma \ref{Tightness}.} According to Theorem 15.5 in \cite{Bil99}, the sequence $(\Pp_a^{(t)}(\cdot))$ is tight if and only if these two conditions hold:\\
\begin{description}
\item[(i)] {\it For each positive $\eta$, there exists an $A$ such that 
\begin{equation}
\Pp_a\left( \frac{|{Z}_{0}|}{\sqrt{t}} \geq A   \big| \inf_{s\leq t} {Z}_s\geq 0 \right)\leq \eta,\qquad t\geq 0.
\end{equation}
}

\item[(ii)]{\it  For each positive $\epsilon$ and $\eta$, there exist $\delta\in (0,1)$, and $t_0>0$ such that 
\begin{equation}
\label{tightness2}\Pp_a\left(  \sup_{s_1,s_2\in [0,1],\, |s_1-s_2|\leq \delta}{t^{-1/2}}|{Z}_{s_1 t}-{Z}_{s_2 t}| \geq \epsilon\big|  \inf_{s\leq t} {Z}_s\geq 0   \right)\leq \eta,\qquad \forall t\geq t_0.
\end{equation}}
\end{description}
Assertion (i) is trivially satisfied. Let us prove (ii), here we follow an idea from  \cite{Bel72}. By (\ref{2.5}), we know that there exists $c>0$ such that
\begin{eqnarray*}
\lefteqn{ \Pp_a\left(  \sup_{s_1,s_2\in [0,1],\, |s_1-s_2|\leq \delta}{t^{-1/2}}|{Z}_{s_1 t}-{Z}_{s_2 t}| \geq \epsilon\big|  \inf_{s\leq t} {Z}_s\geq 0   \right) }\\
&& \leq c\sqrt{t}\Pp_a\left(  \sup_{s_1,s_2\in [0,1],\, |s_1-s_2|\leq \delta}{t^{-1/2}}|{Z}_{s_1 t}-{Z}_{s_2 t}| \geq \epsilon,\, \inf_{s\leq t} {Z}_s\geq 0  \right)  .
\end{eqnarray*}
Let $\epsilon>0$.  For every $\tau \in (0,1)$ and $\delta<\tau$, 
\begin{eqnarray*}
\lefteqn{\Pp_a\left(  \left.\sup_{s_1,s_2\in [0,1],\, |s_1-s_2|\leq \delta}{t^{-1/2}}|{Z}_{s_1 t}-{Z}_{s_2 t}| \geq \epsilon\right|  \inf_{s\leq t} {Z}_s\geq 0   \right)}
\\
&&\leq c\sqrt{t} \Bigg[\Pp_a\left(  \sup_{ s_1,s_2\in [0,1]^2 |s_1-s_2|\leq \delta}{t^{-1/2}}|{Z}_{s_1 t}-{Z}_{s_2 t}| \geq \epsilon,\, \sup_{s\leq \tau}{Z}_{st}< \epsilon \sqrt{t}   ,\,  \inf_{s\leq t} {Z}_s\geq 0  \right) \\
&& + \Pp_a\left( \sup_{s\leq \tau}{Z}_{st}\geq \epsilon \sqrt{t},\,  \inf_{s\leq t} {Z}_s\geq 0  \right)            \Bigg].
\end{eqnarray*}
For the first term we claim that
\begin{equation}
\lim_{\tau\to 0} \limsup_{t\to\infty} \sqrt{t} \Pp_a\left( \sup_{s\leq \tau}{Z}_{st}\geq \epsilon \sqrt{t},\,  \inf_{s\leq t} {Z}_s\geq 0  \right) =0.
\end{equation}
Let $T:=\inf\{s\geq 0,\, {Z}_s\geq \epsilon \sqrt{t}\}.$ 
\begin{eqnarray*}
\Pp_a\left( \sup_{s\leq t\tau}{Z}_s\geq \epsilon \sqrt{t}  ,\, \inf_{s\leq t}{Z}_s\geq 0 \right)= \Pp_a\left(\sup_{s\leq t\tau} {Z}_s\geq \epsilon \sqrt{t},\, \inf_{s\leq t\tau}{Z}_s \geq 0,\, T> t\tau -M\right)
\\
+\Pp_a\left( \sup_{s\leq t\tau} {Z}_s\geq \epsilon \sqrt{t},\inf_{s\leq t\tau}{Z}_s\geq 0,\, T\leq t\tau-M \right).
\end{eqnarray*}
By decomposing on the value of ${Z}_{t\tau-M}$ we have
\begin{eqnarray*}
\lefteqn{\Pp_a\left( \sup_{s\leq t\tau} {Z}_s\geq \sqrt{t}\epsilon,\, \inf_{s\leq t\tau}{Z}_s\geq 0,\, T>t\tau-M   \right)}
\\
&&\leq \Pp_a\left( \frac{{Z}_{t\tau-M}}{\sqrt{t}}\geq \epsilon/2,\, \inf_{s\leq t\tau-M}{Z}_s\geq 0\right)\\
&&+ \Pp_a\left( \sup_{s\in [0,M]}\frac{{Z}_s}{\sqrt{t}}\geq \epsilon/2  \right)\Pp_a\left(\inf_{s\leq t\tau-M}{Z}_s\geq  0 \right)  .
\end{eqnarray*}
Otherwise oberve that ${Z}_T\geq \epsilon \sqrt{t}$ and $t\tau -T\geq M$ imply
\begin{eqnarray*}
&&\Pp_{{Z}_T}\left( {Z}_{t\tau-T}\geq {Z}_T,\, \inf_{s\leq t\tau -T}{Z}_s\geq 0  \right)\leq \Pp\left( {Z}_u\geq 0,\, \inf_{s\leq u} {Z}_s\geq -\epsilon \sqrt{t} \right)_{\big| u=t\tau -T\geq M}
\\
&&\leq \Pp\left( {Z}_u\geq 0,\, \inf_{s\leq u}{Z}_s \geq -\epsilon \sqrt{t\tau}/ \sqrt{\tau}  \right)_{\big| u=t\tau -T\geq M}
\\
&&\leq \Pp\left({Z}_u\geq 0,\, \inf_{s\leq u} {Z}_s \geq -\epsilon \frac{\sqrt{u}}{\sqrt{\tau}}  \right).
\end{eqnarray*}
So when $t$ is large enough, by the Proposition \ref{Donsker22} we have
\begin{eqnarray*}
&&\Pp_a\left( \sup_{s\leq t\tau}{Z}_s \geq\sqrt{t}\epsilon,\, \inf_{s\leq t\tau}{Z}_s\geq 0,\, T<t\tau-M \right)
\\
&&=\mathbb{E}_a\left(  1_{\{ T<t\tau-M \}} 1_{\{ \inf_{s\leq t\tau} {Z}_s\geq 0  \}}  1_{\{ \sup_{s\leq t\tau}{Z}_s \geq \sqrt{t} \epsilon  \}}  \frac{\Pp_{{Z}_T}\left(    {Z}_{t\tau-T}\geq {Z}_T,\, \inf_{s\leq t\tau -T}{Z}_s\geq 0   \right)}{\frac{1}{2}\lambda(\tau,\epsilon)}  \right).
\end{eqnarray*}
with $\lambda(\delta,\epsilon):= \lim_{t\to\infty} \Pp\left( {Z}_t\geq 0,\, \inf_{s\leq t} {Z}_s\geq -\frac{\epsilon \sqrt{t}}{\sqrt{\delta}}  \right)$. Now by using the Markov property at time $T$ we get
\begin{eqnarray*}
\lefteqn{\Pp_a\left(  \sup_{s\leq t\tau} {Z}_s \geq \sqrt{t}\epsilon ,\, \inf_{s\leq t\tau} {Z}_s\geq 0,\, T<t\tau-M \right)}
\\
&&\leq \frac{2}{\lambda(\tau,\epsilon)} \mathbb{E}_a\left( 1_{\{   T<t\tau-M,\, \inf_{s\leq t\tau-M}{Z}_s\geq 0  \}}\Pp_{{Z}_T}\left(  {Z}_{t\tau-T}\geq {Z}_T,\, \inf_{s\leq t\tau -T}{Z}_s\geq 0  \right)  \right)
\\
&&= \frac{2}{\lambda(\tau,\epsilon)} \Pp_a\left( T<t\tau-M,\, {Z}_{t\tau} \geq {Z}_T,\, \inf_{s\leq t\tau}{Z}_s\geq 0  \right)
\\
&&\leq \frac{2}{\lambda{(\tau,\epsilon)}} \Pp_a\left( {Z}_{t\tau}\geq \epsilon\sqrt{t},\, \inf_{s\leq t\tau}{Z}_s\geq 0\right) .
\end{eqnarray*}
Now let us bound the first term. On $\{ \sup_{s\leq \tau}{Z}_{st}< \epsilon \sqrt{t},\, \inf_{s\leq t} {Z}_s\geq 0 \}$ it is clear (note that ${Z}_s\geq 0$ implies $|{Z}_{s_1t}-{Z}_{s_2t}|\leq \max({Z}_{s_1t},{Z}_{s_2t})$) that 
$$\sup_{ s_1,s_2\in [0,\tau]^2,\,  |s_1-s_2|\leq \delta}{t^{-1/2}}|{Z}_{s_1 t}-{Z}_{s_2 t}| \leq \epsilon. $$
Then it follows that
\begin{eqnarray*}
&&\Pp_a\left(  \sup_{ s_1,s_2\in [0,1]^2, |s_1-s_2|\leq \delta}{t^{-1/2}}|{Z}_{s_1 t}-{Z}_{s_2 t}| \geq \epsilon,\, \sup_{s\leq \tau}{Z}_{st}< \epsilon \sqrt{t}   ,\,  \inf_{s\leq t} {Z}_s\geq 0  \right)
\\
&&\leq \Pp_a\left(    \sup_{ s_1,s_2\in [\tau-\delta,1]^2, |s_1-s_2|\leq \delta}{t^{-1/2}}|{Z}_{s_1 t}-{Z}_{s_2 t}| \geq \epsilon  ,\,  \inf_{s\leq t(\tau-\delta)} {Z}_s\geq 0  \right).
\end{eqnarray*}
By the Markov property at time $t(\tau-\delta)$ and (\ref{2.5}), we have
\begin{eqnarray}
 &&\frac{\sqrt{t}}{a} \Pp_a\left(  \sup_{ s_1,s_2\in [0,1]^2, |s_1-s_2|\leq \delta}{t^{-1/2}}|{Z}_{s_1 t}-{Z}_{s_2 t}| \geq \epsilon,\, \sup_{s\leq \tau}{Z}_{st}< \epsilon \sqrt{t}   ,\,  \inf_{s\leq t} {Z}_s\geq 0  \right)
 \notag\\
&& \leq \frac{\sqrt{t}}{a} \Pp_a\left(  \inf_{s\leq t(\tau-\delta)} {Z}_s\geq 0  \right) \Pp\left(  \sup_{ s_1,s_2\in [\tau-\delta,1]^2, |s_1-s_2|\leq \delta}{t^{-1/2}}|{Z}_{s_1 t}-{Z}_{s_2 t}| \geq \epsilon   \right)
\notag \\
&& \leq \frac{1}{\sqrt{\tau-\delta}}\Pp\left( \sup_{ s_1,s_2\in [0,1]^2, |s_1-s_2|\leq \delta}{t^{-1/2}}|{Z}_{s_1 t}-{Z}_{s_2 t}| \geq \epsilon \right).
\label{tozero}
\end{eqnarray}
By Proposition \ref{Donsker22}, the sequence ${t^{-1/2}}({Z}_{st})_{s\in [0,1]}$ is tight. We then have that
\begin{equation}
\lim_{\delta\to0} \limsup_{t\to \infty} \Pp\left( \sup_{ s_1,s_2\in [0,1]^2 |s_1-s_2|\leq \delta}{t^{-1/2}}|{Z}_{s_1 t}-{Z}_{s_2 t}| \geq \epsilon \right)=0.
\end{equation} 
Applying this on the right-hand side of (\ref{tozero}) concludes the proof.
\hfill$\Box$
\begin{Lemma}
\label{finitedim}
Let $a>0$. For any $(s_1<...<s_d)\in [0,1]^d$, we have 
\begin{equation}
\label{eqfinitedim}
\Pp_a^{(t)}(\pi^{-1}_{s_1,...,s_k})\Rightarrow \mathbb{W}^+(\pi^{-1}_{s_1,...,s_d}).
\end{equation}
\end{Lemma}
\noindent{\it Proof of Lemma \ref{finitedim}.} The intuitive idea behind this proof is simple. For any $\kappa>0$, the process $({S}^{(\kappa)}_i)_{i\geq 0}:= (({Z}_{i\frac{1}{\kappa}})_{i\geq 0})_{\kappa\in \N}$ is a standard random walk and then we can apply the result of  \cite{Igl74} to this process.  Moreover, for $\kappa>0$ large enough $ ({S}^{(\kappa)}_i)_{i\geq 0}$ is a good approximation of $({Z}_s)_{s\geq 0}$.\\

 To lighten the proof, we assume that $t=n\in \N$, the general case is similar. Let $d\in \N^*$,  let $F$ be a continuous function from $\mathbb{R}^d$ to $\mathbb{R}$ and let $(s_1,...,s_d) \in [0,1]\cap\mathbb{Q}$. Let $s^*:= \sup\{ s>0,\, \forall j\in [|1,d|],\, \exists p\in \N,\, s_j=ps^*,\, \frac{1}{s^*}\in \N \}$. Observe that for any $\kappa \in \N^*$ the process $({S}^{(\frac{s^*}{\kappa} )}_i)_{\in \N}:= ({Z}_{\frac{s^*}{\kappa}i})_{i\in \N}$ is a standard random walk whose the law of the step is ${Z}_{\frac{s^*}{\kappa}}$, in particular $ \mathbb{E}(({S}_1^{(\frac{s^*}{\kappa})}))=0$ and $\mathbb{E}(({S}_1^{(\frac{s^*}{\kappa})})^2)= \frac{s^*}{\kappa} \sigma^2 $. Then by Lemma 5.1 in \cite{Mad12}, for any $a>0$, 
\begin{eqnarray*}
\lim_{n\to \infty} \Ee_a(F( \frac{1}{\sqrt{ n}} {S}^{(\kappa)}_{s_1\frac{\kappa }{s^*}n},...,\frac{1}{\sqrt{n}} {S}^{(\kappa)}_{s_d\frac{ \kappa}{s^*} n}) \big| \min_{i\in [1,\frac{\kappa}{s^*}n] } {S}^{(\kappa)}_i \geq 0 ) &=& \Ee(F(\sigma m_{s_1},..., \sigma m_{s_d})  ).
\end{eqnarray*}
Observe that
$$  \Ee_a(F( \frac{{S}^{(\kappa)}_{s_1\frac{\kappa }{s^*}n}}{\sqrt{ n}} ,...,\frac{ {S}^{(\kappa)}_{s_d\frac{ \kappa}{s^*} n}}{\sqrt{n}} ) \big| \min_{i\in [1,\frac{\kappa}{s^*}n] } {S}^{(\kappa)}_i \geq 0 ) = \Ee_a(F( \frac{{Z}_{s_1 n}}{\sqrt{n}} ,..., \frac{{Z}_{s_d n}}{\sqrt{n}}) \big| \min_{i\in [1,\frac{\kappa}{s^*}n]} {Z}_{\frac{s^*}{\kappa}i} \geq 0 ).  $$
By the triangle inequality, we have
\begin{eqnarray*}
&&\left|  \Ee_a(F( \frac{ {Z}_{s_1 n} }{\sqrt{t}} ,..., \frac{{Z}_{s_d n}  }{\sqrt{n}}) \big| \min_{i\in [1,\frac{\kappa}{s^*}n]} {Z}_{\frac{s^*}{\kappa}i} \geq 0 )-\Ee_a(F( \frac{{Z}_{s_1 n}}{\sqrt{n}} ,..., \frac{{Z}_{s_d n}}{\sqrt{n}}) \big| \inf_{s\leq n}{Z}_s \geq 0 )\right|
\\
&&\leq ||F|| \Ee_a\left(  \left| \frac{1_{\{ \min_{i\in [1,\frac{\kappa}{s^*}n]} {Z}_{\frac{s^*}{\kappa}i} \geq 0  \}}}{\Pp_a(\min_{i\in [1,\frac{\kappa}{s^*}n]} {Z}_{\frac{s^*}{\kappa}i} \geq 0 ) } -  \frac{1_{\{   \inf_{s\leq n}{Z}_s \geq 0 \}}}{\Pp_a(  \inf_{s\leq n}{Z}_s \geq 0) } \right|        \right)
\\
&&\leq c \frac{\Pp_a\left(\min_{i\in [1,\frac{\kappa}{s^*}n]} {Z}_{\frac{s^*}{\kappa}i} \geq 0   \right)- \Pp_a\left(  \inf_{s\leq n}{Z}_s \geq 0  \right)}{\Pp_a\left(\min_{i\in [1,\frac{\kappa}{s^*}n]} {Z}_{\frac{s^*}{\kappa}i} \geq 0   \right)} + c \left[ 1 -\frac{\Pp_a\left(   \inf_{s\leq n}{Z}_s \geq 0\right)}{\Pp_a\left( \min_{i\in [1,\frac{\kappa}{s^*}n]} {Z}_{\frac{s^*}{\kappa}i} \geq 0   \right)}\right]
\\
&&\leq 2c\left[  1 -\frac{\Pp_a\left(   \inf_{s\leq n}{Z}_s \geq 0\right)}{\Pp_a\left( \min_{i\in [1,\frac{\kappa}{s^*}n]} {Z}_{\frac{s^*}{\kappa}i} \geq 0   \right)} \right].
\end{eqnarray*}
Moreover we claim that
\begin{Lemma}
\label{leclaim1}
For any $a,s^*>0$, 
\begin{equation}
\label{claim1} \lim_{\kappa\to\infty}\lim_{n\to\infty} \frac{\Pp_a(\min_{i\in [1,\frac{\kappa}{s^*}n]} {Z}_{\frac{s^*}{\kappa}i} \geq 0  )}{\Pp_a( \inf_{s\leq n}{Z}_s\geq 0 )}=1.
\end{equation}
\end{Lemma}
The proof of Lemma \ref{leclaim1} is postponed at the end of this  proof.  We deduce that for any $(s_1,...,s_d)\in \mathbb{Q}\cap [0,1]$, when $n$ then $\kappa$ go to infinity 
\begin{equation}
\label{rationel}
  \lim_{n\to\infty} \left| \Ee_a(F( \frac{{Z}_{s_1 n}}{\sqrt{n}} ,..., \frac{{Z}_{s_d n}}{\sqrt{n}}) \big| \inf_{s\leq n} {Z}_s \geq 0 )-  \Ee(F(\sigma m_{s_1},..., \sigma m_{s_d})  )\right| =0.
\end{equation}
Now let $(s_1,...,s_d)\in [0,1]^d$ and $\epsilon>0$. By Proposition (\ref{Donsker22}) for  L\'evy process, we can affirm that
\begin{equation}
\lim_{A\to\infty} \sup_{t>0} \Pp_a\left( \sup_{s\leq 1} |{Z}_{st}|\geq A\sqrt{t}   \right)=0.
\end{equation}
Thus we can assume that $F$ is uniformly continuous. Let $\delta_1>0$ such that for any $((u_1,...,u_d), (v_1,...,v_d)) \in (\r^d)^2$, 
\begin{eqnarray}
\sum_{i=1}^d |u_i-v_i|\leq \delta_1 \Longrightarrow \big| F(u_1,...,u_d)-F(v_1,...v_d)|\leq \epsilon.
\end{eqnarray}
Let $(s_1^{(\delta)},...,s_d^{(\delta)})\in ([0,1]\cap \mathbb Q)^d$ such that
$\sum_{i=1}^d |s_i-s_i^{(\delta)}|\leq \delta_2$. Observe that
\begin{eqnarray*}
&&\left|\Ee_a(F( \frac{{Z}_{s_1 n}}{\sqrt{n}} ,..., \frac{{Z}_{s_d n}}{\sqrt{n}}) \big| \inf_{s\leq n} {Z}_s \geq 0 )- \Ee_a(F( \frac{{Z}_{s_1^{(\delta)} n}}{\sqrt{n}} ,..., \frac{{Z}_{s_d^{(\delta)} n}}{\sqrt{n}}) \big| \inf_{s\leq n} {Z}_s \geq 0 )\right|
\\
&&\leq \epsilon + \Pp_a\left( \sup_{i\in [|1,d|]}|{Z}_{s_i t}-{Z}_{s_i^{(\delta)}t}|\geq \sqrt{t} \delta_1\big|  \inf_{s\leq n} {Z}_s \geq 0  \right)
\\
&&\leq \epsilon + \Pp_a\left(  \sup_{s_1,s_2\in [0,1],\, |s_1-s_2|\leq \delta_2} |{Z}_{s_1 t}-{Z}_{s_2 t}| \geq \sqrt{t}  \delta_1 \big|  \inf_{s\leq t} {Z}_s\geq 0   \right).
\end{eqnarray*}
By applying (\ref{tightness2}), we can choose  $\delta_2$ small and $n$ so that
\begin{equation}
\left|\Ee_a(F( \frac{{Z}_{s_1 n}}{\sqrt{n}} ,..., \frac{{Z}_{s_d n}}{\sqrt{n}}) \big| \inf_{s\leq n} {Z}_s \geq 0 )- \Ee_a(F( \frac{{Z}_{s_1^{(\delta)} n}}{\sqrt{n}} ,..., \frac{{Z}_{s_d^{(\delta)} n}}{\sqrt{n}}) \big| \inf_{s\leq n} {Z}_s \geq 0 )\right|\leq 2\epsilon.
\end{equation}
Similarly for $\delta_2>0$ small enough we have
\begin{equation}
\left|\Ee_a(F(\sigma m_{s_1},...,\sigma m_{s_d} ))- \Ee_a(F((\sigma m_{s_1^{(\delta_2)}},...,\sigma m_{s_d^{(\delta_2)}}  )\right|\leq 2\epsilon.
\end{equation}
So by the triangle inequality, and (\ref{rationel}), we deduce that, for any $\epsilon>0$, there exists $n$ large enough such that
\begin{equation}
\left|\Ee_a(F( \frac{1}{\sqrt{n}} {Z}_{s_1 n},..., \frac{1}{\sqrt{n}}{Z}_{s_d n}) \big| \inf_{s\leq n} {Z}_s \geq 0 ) -  \Ee_a(F(\sigma m_{s_1},...,\sigma m_{s_d} ))\right|\leq \epsilon.
\end{equation}
This concludes the proof of Lemma \eqref{finitedim}. \hfill$\Box$
\\

\noindent{\it Proof of Lemma \ref{leclaim1}.} For any $A>0$, observe that
\begin{eqnarray}
\nonumber &&\Pp_a(\min_{i\in [|1,n k|]} {Z}_{\frac{i}{k}}\geq 0)-\Pp_a(\inf_{s\leq n} {Z}_s\geq 0)
\\
\nonumber &&=  \Pp_a\left( \min_{i\in [|1,n k|]} {Z}_{\frac{i}{k}}\geq 0,\,\inf_{s\leq n} {Z}_s< 0 \right)
\\
\label{claim1term}&&\leq \Pp_a\left(   \min_{i\in [|1,n k|]} {Z}_{\frac{i}{k}}\geq 0,\,\inf_{s\in [A,n]} {Z}_s< 0   \right) +  \Pp\left( \min_{i\in [|1,n k|]} {Z}_{\frac{i}{k}}\geq 0,\,\inf_{s\leq A} {Z}_s< 0 \right).
\end{eqnarray}
Let us study the first term on the right-hand side above. 
\begin{eqnarray*}
&&\Pp_a\left(   \min_{i\in [|1,n k|]} {Z}_{\frac{i}{k}}\geq 0,\,\inf_{s\in [A,n]} {Z}_s< 0   \right) \leq \Pp_a\left(   \min_{i\in [|1,n |]} {Z}_{i}\geq 0,\,\inf_{s\in [A,n]} {Z}_s< 0   \right) 
\\
&&\leq  \sum_{i=A}^n\Pp_a( \min_{i\in [|1,n|]}{Z}_i \geq 0,\, \inf_{s\in [i,i+1]}{Z}_s<0)
\\
&&\leq \sum_{i=A}^n \sum_{k,p=0}^{+\infty} \Pp_a\left( \min_{i\in [|1,n|]} {Z}_i\geq 0,\, {Z}_i\in [k,k+1],\, {Z}_{i+1}\in [p,p+1],\, \inf_{s\in [i,i+1]}({Z}_s-{Z}_i)\leq -k \right).
\end{eqnarray*}
By the markov property at time $i$ then $i+1$ we get
\begin{eqnarray*}
&&\Pp_a\left(   \min_{i\in [|1,n k|]} {Z}_{\frac{i}{k}}\geq 0,\,\inf_{s\in [A,n]} {Z}_s< 0   \right)
\\
&&\leq \sum_{i=A}^n \sum_{k,p=0}^{+\infty} \Pp_a\left( \min_{j\in [|1,i|]} {Z}_i\geq 0,\, {Z}_i\in [k,k+1]\right)\\
&&\hspace{3cm}\times \Pp\left(|{Z}_{1}| \geq |k-p|-2 \inf_{s\in [0,1]}{Z}_s\leq -k \right)\Pp_p(  \min_{j\in [|1,n-i|]} {Z}_{j}\geq 0 ).
\end{eqnarray*}
By using (\ref{2.5}), (\ref{2.6}) and exponential moments of ${Z}$, we get, for an appropriate choice of $\epsilon_0>0$,
\begin{eqnarray*}
\Pp_a\left(   \min_{i\in [|1,n k|]} {Z}_{\frac{i}{k}}\geq 0,\,\inf_{s\in [A,n]} {Z}_s< 0   \right)   &\leq& \sum_{i=A}^n \sum_{k,p=0}^{+\infty} \frac{(1+a)(1+k)}{i^\frac{3}{2}} \ee^{-\epsilon_0 k} \ee^{-\epsilon_0 |p-k|} \frac{(1+p)}{\sqrt{n-i}}
\\
&\leq& C(1+a)\sum_{i=A}^n \frac{1}{i^\frac{3}{2} \sqrt{n-i}} \leq \frac{C(1+a)}{A\sqrt{n}}.
\end{eqnarray*}
We fix $A$ large enough such that $\frac{C(1+a)}{A}\leq \epsilon$. Now we treat the second term in (\ref{claim1term}). By the Markov property at time $A$ we get
\begin{eqnarray*}
\Pp\left( \min_{i\in [|1,n k|]} {Z}_{\frac{i}{k}}\geq 0,\,\inf_{s\leq A} {Z}_s< 0 \right)\leq  \sum_{p\geq 0}\Pp\left( \min_{i\in [|1,n k|]} {Z}_{\frac{i}{k}}\geq 0,\,\inf_{s\leq A} {Z}_s< 0,\, {Z}_A\in [p,p+1] \right) .
\end{eqnarray*}
By the Cauchy-Schwarz inequality then  (\ref{2.5}), we get
\begin{eqnarray*}
&& \Pp_a\left( \min_{i\in [|1,n k|]} {Z}_{\frac{i}{k}}\geq 0,\,\inf_{s\leq A} {Z}_s< 0 \right) 
\\
&&\leq \sum_{p\geq 0} \Pp_a\left( \min_{i\in [|1,A k|]} {Z}_{\frac{i}{k}}\geq 0,\,\inf_{s\leq A} {Z}_s< 0,\, {Z}_A\in [p,p+1] \right) \Pp_{p+1}\left( \min_{i\in [|1,nk-Ak|]}{Z}_{\frac{i}{k}} \geq 0 \right)
\\
&&\leq \sum_{p\geq 0} (1+p)\Pp_a\left( {Z}_A\in [p,p+1] \right)^\frac{1}{2}  \frac{c}{\sqrt{n-A}}\Pp\left(   \min_{i\in [|1,A k|]} {Z}_{\frac{i}{k}}\geq 0,\,\inf_{s\leq A} {Z}_s< 0 \right)^\frac{1}{2}
\\
&&\leq \frac{c'}{\sqrt{n-A}} \Pp\left(   \min_{i\in [|1,A k|]} {Z}_{\frac{i}{k}}\geq 0,\,\inf_{s\leq A} {Z}_s< 0 \right)^\frac{1}{2}.
\end{eqnarray*}
Clearly for any fixed $A>0$ $\lim_{k\to\infty}\Pp\left(   \min_{i\in [|1,A k|]} {Z}_{\frac{i}{k}}\geq 0,\,\inf_{s\leq A} {Z}_s< 0 \right)=0$. It achieves the proof of (\ref{claim1}).
\hfill$\Box$

%
%

 \renewcommand{\theequation}{B.\arabic{equation}}
  \renewcommand{\thetheorem}{B.\arabic{theorem}}

  \setcounter{equation}{0}  
    \setcounter{theorem}{0}  

 \subsection{Two inequalities}
 \begin{Lemma}
 \label{2ineq}
 There exist constant $c_1,c_2,c_3$ such that for any $a\geq 0, u\geq v\geq 0$ and $t\geq 1$,
 \begin{eqnarray}
 \label{2.5}  \Pp_a\left( \inf_{s\leq t} {Z}_s \geq 0\right)& \leq &c_1 (1+a) t^{-\frac{1}{2}},
 \\
 \label{2.6} \Pp_a\left( {Z}_t\in [u,v],\, \inf_{s\leq t} {Z}_s \geq 0 \right)   &\leq& c_2(1+a) (1+v-u) (1+v) t^{-\frac{3}{2}},
 \\
 \label{toadd} \mathbb E_a\left(  Z_t \1_{\{ \inf_{s\leq t} Z_s\geq 0 \}} \right)&\leq &c_3 a.
 \end{eqnarray}
 \end{Lemma}
 \noindent{\it Proof of Lemma \ref{2ineq}.} Under the general assumptions we have made on $Z$, the restriction of $({Z}_s)_{s\geq 0}$ to the integers is a centered random walk with $0$ mean and finite variance. Therefore, according to the lemmas in the Section 2 of \cite{Aid11}, (\ref{2.5}) follows. Inequality (\ref{2.6}) follows with the same spirit. It suffices to observe that
 \begin{eqnarray*}
 \Pp_a({Z}_t \in [u,v],\, \inf_{s\leq t} {Z}_s \geq 0) &\leq &\Pp_a(  \min_{i\in [|1,\lfloor t|]} {Z}_i\geq 0,\, {Z}_t\in [a,b]    )
 \\
 &\leq& \sum_{k\geq 0}\Pp_a(  \min_{i\in [|1,\lfloor t|]} {Z}_i,\, {Z}_{\lfloor t\rfloor} \in [k,k+1]    )\max_{y\in [k,k+1]}\Pp( {Z}_{t-\lfloor t\rfloor} \in [u-y,v-y]).
 \end{eqnarray*}
 Using (\ref{chap3Aine}) we get that
 \begin{eqnarray*}
 \Pp_a({Z}_t \in [u,v],\, \inf_{s\leq t} {Z}_s \geq 0)&\leq &   c (1+a) t^{-\frac{3}{2}}  \sum_{k\geq 0} (1+k) \sup_{y\in [k,k+1]}\Pp( {Z}_{t-\lfloor t\rfloor} \in [u-y,v-y]).
 \end{eqnarray*}
 Moreover as $\mathbb{E}(\ee^{\epsilon_0 |{Z}_1})<\infty$, by the Markov inequality, we have that
 \begin{eqnarray*}
  \Pp_a({Z}_t \in [u,v],\, \inf_{s\leq t} {Z}_s \geq 0)&\leq & c (1+u) t^{-\frac{3}{2}} \left\{   \sum_{k= 0}^u (1+k) \ee^{-\epsilon_0(u-k)}  +  \sum_{k= v}^{+\infty} (1+k) \ee^{-\epsilon_0(k-v)} + \sum_{k=u}^v (1+k)       \right\}
 \\
 &\leq &c_2 (1+a) (1+v-u) (1+v) t^{-\frac{3}{2}},
 \end{eqnarray*}
 which ends the proof of \eqref{2.6}. By similar arguments and Lemma 2.3 in \cite{AShi11} one can obtain (\ref{toadd}). The details are left to the reader  
 \hfill$\Box$

  \begin{Lemma}
  \label{PourRW2}
 For any $\theta>0$, there exists a constant $c(\theta )>0$ such that for any $a\geq 0$,
  \begin{eqnarray}
  	\label{eqPourRW2}
 \int_0^\infty    \mathbb{E}_a\left( \1_{\{  \inf_{u\leq s} {Z}_s  \geq 0 \}}  \ee^{-\theta {Z}_s}      \right)  ds \leq c(\theta).
  \end{eqnarray}
  \end{Lemma}
  \noindent{\it Proof of Lemma \ref{PourRW2}.}We observe that
  \begin{eqnarray*}
  \int_0^\infty    \mathbb{E}_a\left( \1_{\{  \inf_{u\leq s} {Z}_s  \geq 0 \}}  \ee^{-\theta {Z}_s}      \right) &\leq  & \int_0^\infty   \mathbb{E}_a\left( \1_{\{  \inf_{i\in [|1,\lfloor s\rfloor |] } {Z}_i  \geq 0 \}}  \ee^{-\theta {Z}_i}\mathbb{E}\left(  \ee^{-\theta {Z}_{s-\lfloor s\rfloor }}  )      \right) \right)ds
  \\
  &\leq & c \sum_{i=0}^{\infty} \mathbb{E}_a\left( \1_{\{  \inf_{i\in [|1,\lfloor s\rfloor |] } {Z}_i  \geq 0 \}}  \ee^{-\theta {Z}_i}  \right),
  \end{eqnarray*}
  and we conclude by Lemma B.2 in \cite{Aid11}. 
  \hfill$\Box$

\bibliographystyle{plain}
\bibliography{bibli}

\end{document}